\documentclass{tac}

\usepackage[square]{natbib} 
\usepackage{mathtools} 
\usepackage{amsfonts} 
\usepackage{enumitem} 

\usepackage[colorlinks=true]{hyperref}
\hypersetup{allcolors=[rgb]{0.1,0.1,0.4}}

\author{Marco Abbadini}

\address{Department of Mathematics {\sl Federigo Enriques}, Universit\`a degli Studi di Milano\\
	via Cesare Saldini 50, 20133 Milano, Italy.}

\title {The dual of compact ordered spaces is a variety}

\copyrightyear{2019}

\keywords{compact ordered space, variety, duality, axiomatizability.}
\amsclass{Primary: 03C05. Secondary: 08A65, 18B30, 18C10, 54A05, 54F05}

\eaddress{marco.abbadini@unimi.it}

\mathrmdef{C}
\mathrmdef{Max}
\mathrmdef{op}
\mathrmdef{ar}
\mathrmdef{ev}
\mathrmdef[dist]{d}
\newcommand{\du}{{\dist^\uparrow}}
\mathopdef{essinf}
\mathopdef{esssup}
\mathopdef{ess}

\mathbfdef{Set}
\mathbfdef{PosComp}
\mathbfdef[KH]{CompHaus}
\newcommand{\MCinfty}{\mathbf{MC}_\mathbf{\infty}}
\mathbfdef{MC}
\mathbfdef[PreT]{PreTop}

\newcommand{\N}{\mathbb{N}}

\mathcaldef[Zm]{Z}
\newcommand{\Zt}{\widetilde{\Zm}}
\mathcaldef{A}

\begin{document}
	\maketitle
	\begin{abstract}
		In a recent paper (2018), D.\ Hofmann, R.\ Neves and P.\ Nora proved that the dual of the category of compact ordered spaces and monotone continuous maps is a quasi-variety---not finitary, but bounded by $\aleph_1$. An open question was: is it also a variety? We show that the answer is affirmative. We describe the variety by means of a set of finitary operations, together with an operation of countably infinite arity, and equational axioms. The dual equivalence is induced by the dualizing object $[0,1]$.
	\end{abstract}

\section{Introduction}\label{s:introduction}

Compact ordered spaces were introduced by L.\ Nachbin, and they are to topology and partial order what compact Hausdorff spaces are to topology. A \emph{compact ordered space} $(X,\leq,\tau)$ consists of a compact space $(X,\tau)$ equipped with a partial order $\leq$ so that the set
$$\{(x,y)\in X\times X\mid x\leq y \}$$
is closed in $X\times X$ with respect to the product topology (see \citep{Nachbin} for a standard reference); we are interested in the category \PosComp of compact ordered spaces with monotone continuous maps. 

The goal of this paper is to establish for \PosComp a result which is known to hold for the category \KH of compact Hausdorff spaces with continuous maps---namely, that the dual category is a variety (not finitary, but bounded by $\aleph_1$). As far as compact Hausdorff spaces are concerned, we recall some historical details: in 1969, Duskin proved that the functor $\hom({-}, [0,1] )\colon {\KH}^\op\rightarrow \Set$ is monadic \citep{Duskin}; Isbell presented a set of primitive operations of ${\KH}^\op $, using finitely many finitary operations, along with an operation of countably infinite arity \citep{Isbell}; finally, Marra and Reggio provided finitely many axioms to axiomatize the variety ${\KH}^\op$ \citep{MarraReggio}. 

These results were a source of motivation for the algebraic study of the dual of \PosComp in  \citep{HofmannShort}: the authors proved that $\PosComp^\op$ is a quasi-variety---not finitary, but bounded by $\aleph_1$---leaving as open the following question.
\begin{quote}
	Is $\PosComp^\op $ also a variety?
\end{quote} Our main result is that the answer is affirmative, as stated in the following theorem.
\theorem [Main result]\label{t:MAIN-PosComp-is-variety}
The dual of $\PosComp$ is equivalent to a variety of algebras.
\endtheorem 
\noindent The proof of our main result is at times inspired by \citep{HofmannShort}, but does not depend on their results. In this paper, under the term \emph{variety of algebras}, we admit the so called varieties of \emph{infinitary} algebras, whose operations may have infinite arity (see \citep{Slo} for varieties of infinitary algebras). The variety in Theorem~\ref{t:MAIN-PosComp-is-variety} will be denoted by $\MCinfty$, where $\mathbf{M}$ stands for ``monotone'' and $\mathbf{C}$ stands for ``continuous''. We will give a set of primitive operations (of countable arity) and a set of equational axioms for $\MCinfty$. 

The set $[0,1]$ is both a  compact ordered space, with the canonical order and the euclidean topology, and an MC$_\infty$-algebra, in a natural way. In fact, we present the duality between $\PosComp$ and $\MCinfty$ as induced by the dualizing object $[0,1]$. This dual equivalence coincides, essentially, with the duality available in \citep{HofmannShort}: the main difference is that, on the algebraic side, we consider a slightly different set of primitive operations, that facilitates us to state the axioms in an equational form.

We will show $\MCinfty=\mathrm{ISP}([0,1])$, where $\mathrm{P}$ denotes the closure under products, $\mathrm{S}$ denotes the closure under subalgebras, and $\mathrm{I}$ denotes the closure under isomorphisms. Moreover, we will see that the operations of $\MCinfty$, interpreted in $[0,1]$, are precisely the monotone continuous maps from a power of $[0,1]$ to $[0,1]$. In other words, $\MCinfty$ is the category of algebras of the varietal theory (in the sense of \citep{Linton}) whose objects are powers of $[0,1]$ and whose morphisms are the monotone continuous maps.

\subsection{The strategy}
The strategy that we adopt to prove the duality follows the lines of \citep{MarraReggio}, which, to the best of our knowledge, used this strategy for the first time in a similar context---namely, to obtain a finite equational axiomatization of the dual of $\KH$. We prove that the dual of $\PosComp$ is a variety via the following steps.
\begin{enumerate}
	\item We obtain a dual adjunction between the category of preordered topological spaces and a finitary variety $\MC$, to be defined. This dual adjunction is induced by the dualizing object $[0,1]$.
	\item We characterize the objects which are fixed by the adjunction. On the topological side, the fixed objects are precisely the compact ordered spaces. On the algebraic side, the fixed objects are the archimedean Cauchy complete MC-algebras. Hence, a duality is established between the full subcategories of compact ordered spaces and archimedean Cauchy complete MC-algebras. 
	\item We show that the full subcategory of archimedean Cauchy complete MC-algebras is isomorphic to an infinitary variety $\MCinfty$, obtained by adding a term $\delta$ of countably infinite arity to the language of $\MC$, together with some new appropriate equational axioms. The term $\delta$ is intended to map \textit{enough} Cauchy sequences to their limit. The forgetful functor $\MCinfty\rightarrow \MC$, restricted at codomain, gives the desired isomorphism.
\end{enumerate}
To show that every compact ordered space is fixed by the adjunction, we use an analogue of Urysohn's Lemma. To show that every archimedean Cauchy complete MC-algebra is fixed, we use the Subdirect Representation Theorem, which applies since $\MC$ is finitary, and an analogue of Stone-Weierstrass Theorem. The Subdirect Representation Theorem is used to show that every archimedean MC-algebra is mapped injectively by the unit of the adjunction. The analogue of Stone-Weierstrass Theorem is used to show that if the algebra is Cauchy complete, then it is mapped surjectively by the unit of the adjunction.

\paragraph{\textbf{Acknowledgements.}} The author would like to thank his Ph.D.\ advisor Vincenzo Marra for his suggestions. Moreover, the author is deeply grateful to Luca Reggio, who notably improved the proof of Theorem~\ref{t:X(A) is compact}, and helped with many useful comments. Finally, the author expresses his gratitude to the anonymous referee for his or her careful reading and several comments that helped to achieve a better presentation of the results.

\section{The category \PreT of preordered topological spaces}\label{s:pret}

\definition
	A \emph{preordered topological space} $(X,\leq,\tau)$ consists of a set $X$, a preorder $\leq$ on $X$ and  a topology $\tau$ on $X$.
	
\enddefinition
When no confusion arises, we write $X$ instead of $(X,\leq,\tau)$.
We denote with \PreT the category whose objects are preordered topological spaces and whose morphisms from $A$ to $B$ are the  monotone  (i.e.\ $x\leq y\Rightarrow f(x)\leq f(y)$) continuous functions $f\colon A\rightarrow B$.

It is well known that the forgetful functors from the category of topological spaces and from the category of preordered sets to the category of sets are both topological, in the sense of \citep[Definition 21.1]{Joyofcats}. The forgetful functor $\mathrm{U}_{\PreT}\colon \PreT\rightarrow \Set$ is topological, too; indeed, every $\mathrm{U}_{\PreT}$-structured source $$\left(X\xrightarrow{f_i} \mathrm{U}_{\PreT}(X_i,\leq_i, \tau_i)\right)_{i\in I}$$ admits a unique $\mathrm{U}_{\PreT}$-initial lift 
$$\left((X,\leq,\tau)\xrightarrow{\overline{f}_i}(X_i,\leq_i, \tau_i)\right)_{i\in I},$$
where $\leq$ is defined by 
$$x\leq y\Leftrightarrow \forall i\in I\ f_i(x)\leq f_i(y),$$
and $\tau$ is the topology generated by 
$$\{f_i^{-1}(O_i)\mid i\in I,O_i\in \tau_i \}.$$

It is well known that topological functors lift limits (in particular, products) uniquely \citep[Proposition 21.15]{Joyofcats}. For a family $(X_i,\leq_i)$ of preordered spaces, the \emph{product preorder} on $\prod_{i\in I}X_i$ is the preorder $\leq$ defined by 
$$x\leq y \Longleftrightarrow 	\forall i\in I\ \pi_i(x)\leq_i \pi_i(y),$$
where $\pi_i$ is the projection onto the $i$-th coordinate.
For a family $(X_i,\tau_i)_{i\in I}$ of topological spaces, the \emph{product topology} on $\prod_{i\in I}X_i$ is the topology generated by
$$\{\pi_i^{-1}(O_i)\mid i\in I, O_i\in \tau_i \}.$$
Finally, for a family $(X_i,\tau_i, \leq_i)_{i\in I}$ of preordered topological spaces, the unique $\mathrm{U}_{\PreT}$-initial lift  of the $\mathrm{U}_{\PreT}$-structured source $\left(\prod_{i\in I}X_i\xrightarrow{\pi_i} \mathrm{U}_{\PreT}(X_i,\leq_i, \tau_i)\right)_{i\in I}$ is 
$$\left(\left(\prod_{i\in I} X,\leq,\tau\right)\xrightarrow{\overline{\pi}_i}(X_i,\leq_i, \tau_i)\right)_{i\in I},$$
where $\leq$ is the product preorder, and $\tau$ is the product topology. Moreover, this is a categorical product in $\PreT$. Unless otherwise stated, when referring to the set-theoretic product of preordered topological spaces as a preordered topological space, we implicitly assume that the preorder is the product preorder and the topology is the product topology. Then, if $I$ is a set, and $f\colon [0,1]^I\rightarrow [0,1]$ is a monotone and continuous function, we have the following: for every preordered topological space $X$, $f$ is an internal operation on $\hom_{\PreT}(X,[0,1])$, meaning that, for every $I$-indexed family $(g_i)_{i\in I}$ of monotone continuous functions from $X$ to $[0,1]$, the function $X\rightarrow [0,1]$, $x\mapsto f((g_i(x))_{i\in I})$ is monotone and continuous, as well.

\section{The variety $\MC$}\label{s:V}

We define some operations on $[0,1]$. For $a,b\in [0,1]$, $a\lor b$ and $a\land b$ denote, respectively, the supremum and the infimum of $\{a,b\}$, $a\oplus b\coloneqq \min\{a+b,1 \}$, and $a\odot b\coloneqq\max\{a+b-1,0\}$. Moreover, for each $\lambda\in [0,1]$, the constant symbol $\lambda$ denotes $\lambda$ itself. 
\remark\label{r:all operations in V are mon and cont}
	Each of these operations ($\lor$, $\land$, $\oplus$, $\odot$, and, for every $\lambda\in [0,1]$, the constant function $\lambda$) is monotone and continuous with respect to the product order and product topology.
\endremark 
\noindent Note that we do not consider the function $\lnot\colon [0,1]\rightarrow [0,1]$ $a\mapsto 1-a$, since it is not monotone.

We define a finitary variety $\MC$ of algebras of type $\mathcal{L}=\{\oplus,\odot,\lor,\land,0,1 \}\cup\{\lambda \mid\lambda\in [0,1]\}$. Specifically, an algebra $A$ belongs to $\MC$ (and we say that $A$ is an \emph{MC-algebra}) if it satisfies the following  identities, which, as one may verify, are all satisfied by $[0,1]$.
\begin{enumerate}
	\item	$\langle A,\lor,\land,0,1 \rangle$ is a distributive bounded lattice.
	\begin{enumerate}
		\item $a\lor b=b\lor a$.
		\item $a\land b=b\land a$.
		\item $a\lor(b\lor c)=(a\lor b)\lor c$.
		\item $a\land(b\land c)=(a\land b)\land c$.
		\item $a\lor(a\land b)=a$.
		\item $a\land(a\lor b)=a$.
		\item $a\lor 0=a$.
		\item $a\land 1=a$.
		\item $a\lor(b\land c)=(a\lor b)\land (a\lor c)$.
		\item $a\land (b\lor c)=(a\land b)\lor (a\land c)$.
	\end{enumerate}
	
	\item $\langle A,\oplus,0\rangle$ is a commutative monoid, with absorbing element $1$.		
	\begin{enumerate}
		\item $(a\oplus b)\oplus c=a\oplus(b\oplus c)$.
		\item $(a\oplus b)=(b\oplus a)$.
		\item $a\oplus 0=a $.
		\item $a\oplus1=1$.
	\end{enumerate}
	
	\item $\langle A,\odot,1\rangle$ is a commutative monoid, with absorbing element $0$.		
	\begin{enumerate}
		\item $(a\odot b)\odot c=a\odot(b\odot c)$.
		\item $(a\odot b)=(b\odot a)$.
		\item $a\odot 1=a $.
		\item $a\odot 0=0$.
	\end{enumerate}
	\item $\oplus$ and $\odot$ distribute over $\lor$ and $\land$.
	\begin{enumerate}
		\item $(a\lor b)\oplus c=(a\oplus c)\lor (b\oplus c)$.
		\item $(a\land b)\oplus c=(a\oplus c)\land (b\oplus c).$
		\item $(a\lor b)\odot c=(a\odot c)\lor (b\odot c)$.
		\item $(a\land b)\odot c=(a\odot c)\land (b\odot c).$	
	\end{enumerate}	
	
	\item $ (a\oplus b)\odot c\leq a\oplus (b\odot c)$.
	
	\item For each $\lambda\in [0,1]$, we have the axiom $a\leq (a\odot(1- \lambda))\oplus\lambda$.
	
	\item For each $\lambda\in [0,1]$, we have the axiom $a\geq (a\oplus \lambda)\odot(1- \lambda)$.
	
	\item For every $n,m\in \{0,1,2,\dots\}$, we have the axiom 
	$$a\land (b\oplus \underbrace{(c\odot \lambda)\oplus\dots\oplus(c\odot \lambda)}_{n\text{ times}})\leq (a\odot \underbrace{(c\oplus \lambda)\odot \dots\odot(c\oplus \lambda)}_{m\text{ times}})\lor b.$$
	
	\item \label{i:lor}For $\alpha,\beta,\gamma\in [0,1]$ such that $\alpha\lor \beta=\gamma$ in $[0,1]$, we have the axiom $\alpha\lor \beta=\gamma$.	
	\item \label{i:land}For $\alpha,\beta,\gamma\in [0,1]$ such that $\alpha\land \beta=\gamma$ in $[0,1]$, we have the axiom $\alpha\land \beta=\gamma$.
	\item \label{i:oplus}For $\alpha,\beta,\gamma\in [0,1]$ such that $\alpha\oplus \beta=\gamma$ in $[0,1]$, we have the axiom $\alpha\oplus \beta=\gamma$.
	\item \label{i:odot}For $\alpha,\beta,\gamma\in [0,1]$ such that $\alpha\odot \beta=\gamma$ in $[0,1]$, we have the axiom $\alpha\odot \beta=\gamma$.	
\end{enumerate}
For $\lambda\in [0,1]$, we write $x\ominus \lambda$ for $x\odot (1-\lambda)$. In $[0,1]$, $x\ominus \lambda=\max\{x-\lambda, 0\}$. We remark that we allow the notation $x\ominus \lambda$ only when $\lambda$ is a constant symbol in $[0,1]$.

\section{The dual adjunction between $\PreT$ and $\MC$}

Let $X$ be a preordered topological space. We set
$$
\mld \C(X)&\coloneqq\hom_{\PreT}(X,[0,1])=\\
=\{f\colon X\rightarrow [0,1]\mid f\text{ is monotone and continuous} \}.
$$
Since, by Remark~\ref{r:all operations in V are mon and cont}, the interpretation in $[0,1]$ of every MC-operation is monotone and continuous, $\C(X)$ is an MC-algebra with pointwise defined operations. For each $x\in X$, we set
\begin{align*}
\begin{split}
\ev_x\colon \C(X)&\longrightarrow [0,1]\\
a&\longmapsto a(x).
\end{split}
\end{align*}

Let $A\in\MC$. Set $\Max(A)\coloneqq \hom_{\MC}(A,[0,1])$. The motivation for this name stems from the fact that the set of morphisms from an MC-algebra $A$ to $[0,1]$ is in bijection with the set of maximal congruences on $A$; this follows from the fact that $[0,1]$ is the only simple algebra, as will be proved in Corollary~\ref{c:simple}.
For each $a\in A$, we set 
\begin{align*}
\begin{split}
\ev_a\colon \Max(A)&\longrightarrow [0,1]\\
x&\longmapsto x(a).
\end{split}
\end{align*}
For all $x,y\in \Max(A)$, set $x\leq y$ if, and only if, for all $a\in A$, $\ev_a(x)\leq \ev_a(y)$, i.e., $x(a)\leq y(a)$. Let $\tau$ be the smallest topology on $\Max(A)$ that contains $\ev_a^{-1}(O)$ (i.e., $\{x\in \Max(A)\mid x(a)\in O \}$) for every $a\in A$ and $O$ open subset of $[0,1]$. 

In \citep[Section 1-C]{Tholen}, some properties are discussed that are sufficient to establish a dual adjunction induced by a dualizing object. These properties are expressed in terms of existence of certain initial lifts, and in our case these properties hold. Indeed, let $\mathrm{U}_{\PreT}\colon \PreT\rightarrow \Set$ and $\mathrm{U}_\MC\colon \MC\rightarrow \Set$ denote the forgetful functors; by the results discussed in section \ref{s:pret}, for every $A\in \MC$,  $\left((\Max(A),\leq, \tau)\xrightarrow{\ev_a} [0,1]\right)_{a\in A}$ is the unique $\mathrm{U}_{\PreT}$-initial lift of the $\mathrm{U}_{\PreT}$-structured source $\left(\Max(A)\xrightarrow{\ev_x} \mathrm{U}_\PreT([0,1])\right)_{x\in X}$. Moreover, since for every preordered topological space $X$ the operations in $\C(X)$ are pointwise defined, we have that $$\left(\C(X)\xrightarrow{\ev_x} [0,1]\right)_{x\in X}$$ (where $\C(X)$ denotes the MC-algebra whose underlying set is $\hom_{\PreT}(X,[0,1])$ and with pointwise defined operations) is the unique $\mathrm{U}_{\MC}$-initial lift of the $\mathrm{U}_{\MC}$-structured source $$\left(\C(X)\xrightarrow{\ev_x} \mathrm{U}_\MC([0,1])\right)_{x\in X}$$ (where $\C(X)$ denotes the set $\hom_{\PreT}(X,[0,1])$).

Therefore, we have a dual adjunction between $\PreT$ and $\MC$ induced by the dualizing object $[0,1]$, that we now make explicit. In accordance with \citep{Tholen}, we find a more natural choice to use \emph{contravariant} functors between $\PosComp$ and $\MC$ rather than \emph{covariant} ones between $\PosComp^\op$ and $\MC$, or between $\PosComp$ and $\MC^\op$. This choice seems to us more natural in the context of dual adjunctions induced by a dualizing object, because it respects the symmetry between the two involved categories ($\PosComp$ and $\MC$, in our case). Since we are considering contravariant functors, we end up using \emph{two units}, rather than a unit and a counit.

The assignment $\C$ on the objects may be extended on arrows so that $\C$ becomes a contravariant functor: for a morphism $g\colon X\rightarrow Y$ in $\PreT$, we set
\begin{align*}
\C(g)\colon \C(Y)&\longrightarrow \C(X)\\
a&\longmapsto a\circ g.
\end{align*}
Analogously, the assignment $\Max$ on the objects may be extended on arrows so that $\Max$ becomes a contravariant functor: for a morphism $f\colon A\rightarrow B$ in $\MC$, we set
\begin{align*}
\Max(f)\colon \Max(B)&\longrightarrow \Max(A)\\
x&\longmapsto x\circ f.
\end{align*}
The adjunction is given as follows. Let $X\in \PreT$ and $A\in \MC$. To each morphism $g\colon X\rightarrow \Max(A)$ in $\PreT$ we associate the following morphism in $\MC$: 
\begin{align*}
\hat{g}\colon A&\longrightarrow \C(X)\\
a&\longmapsto \ev_a\circ g;\ x\mapsto (g(x))(a).
\end{align*}
To each morphism $f\colon A\rightarrow \C(X)$ in $\MC$ we associate the following morphism in $\PreT$:
\begin{align*}
\check{f}\colon X&\longrightarrow \Max(A)\\
x&\longmapsto \ev_x\circ f;\ a\mapsto (f(a))(x).
\end{align*}
For $X\in \PreT$, the unit at $X$ is
\begin{align*}
\eta_X\colon X&\longrightarrow \Max\C(X)\\
x&\longmapsto (\ev_x \colon \C(X)\rightarrow [0,1]; a\mapsto a(x)).
\end{align*}
For $A\in \MC$, the unit at $A$ is
\begin{align*}
\varepsilon_A\colon A&\longrightarrow \C\Max(A)\\
a&\longmapsto (\ev_a \colon \Max(A)\rightarrow [0,1]; x\mapsto x(a)).
\end{align*}

\section{Fixed objects on the geometrical side}
Let us recall the definition of compact ordered space.
\definition
	A \emph{compact ordered space} $(X,\leq,\tau)$ consists of a compact space $(X,\tau)$ equipped with a partial order $\leq$ so that the set
	$$\{(x,y)\in X\times X\mid x\leq y \}$$
	is closed in $X\times X$ with respect to the product topology.
\enddefinition 
 A standard reference is \citep{Nachbin}. We recall that every compact ordered space is Hausdorff \citep[Proposition 2, Chapter 1, p.\ 27]{Nachbin}. 
 We denote with \PosComp the category of compact ordered spaces with monotone continuous maps.

The goal of this section is to prove the following.
\theorem \label{t:characterization fixed geometric}
	Let $X$ be a preordered topological space. The following conditions are equivalent.
	\begin{enumerate}
		\item\label{i:unit-iso} The unit $\eta_X\colon X\rightarrow \Max\C(X)$ is an isomorphism.
		\item\label{i:in-the-image} There exists an MC-algebra $A$ such that $X$ and $\Max(A)$ are isomorphic preordered topological spaces.
		\item\label{i:comp-pospace} $X$ is a compact ordered space.
	\end{enumerate}
\endtheorem
\remark\label{r:specialization}
	The implication [\eqref{i:unit-iso}$\Rightarrow$\eqref{i:in-the-image}] in Theorem~\ref{t:characterization fixed geometric} is immediate: take $A= \C(X)$.
\endremark 
\subsection{$\Max(A)$ is a compact ordered space}
In this subsection, we prove the implication [\eqref{i:in-the-image}$\Rightarrow$\eqref{i:comp-pospace}] of Theorem~\ref{t:characterization fixed geometric}, i.e., for every $A\in \MC$, $\Max(A)$ is a compact ordered space. We need the following lemmas and remarks.

\lemma\label{l:PosComp is hereditary}
	Let $X$ be a compact ordered space, and let $Y$ be a closed subset of $X$. Then $Y$, equipped with the topology and the order induced by $X$, is a compact ordered space.
\endlemma
\proof 
	Since $Y$ is a closed subspace of a compact space, $Y$ is compact. Clearly, the partial order on $X$ induces a partial order on $Y$. The product topology on $Y\times Y$ coincides with the subspace topology on $Y\times Y$ as subspace of $X\times X$. Since $\{(x,y)\in X\times X\mid x\leq y \}$ is a closed subset of $X\times X$, $\{(x,y)\in Y\times Y\mid x\leq y \}=\{(x,y)\in X\times X\mid x\leq y \}\cap Y\times Y$ is closed in $Y\times Y$.
\endproof 

\lemma
	Let $(X_i)_{i\in I}$ be a family of compact ordered spaces. Then, $\prod_{i\in I}X_i$, equipped with the product topology and product order, is a compact ordered space.
\endlemma
\proof 
	By Tychonoff's theorem, $\prod_{i\in I}X_i$ is compact. Let us consider the bijection $\phi\colon (\prod_{i\in I}X_i)\times (\prod_{i\in I}X_i)\rightarrow \prod_{i\in I}(X_i\times X_i)$; $((a_i)_{i\in I},(b_i)_{i\in I})\mapsto (a_i,b_i)_{i\in I}$. The function $\phi$ is a homeomorphism, and the image under $\phi$ of the set 
	$$\left\{\left(\left(x_i\right)_{i\in I},\left(y_i\right)_{i\in I} \right)\in  \left(\prod_{i\in I}X_i\right)\times \left(\prod_{i\in I}X_i\right)\mid \left(x_i\right)_{i\in I}\leq(y_i)_{i\in I} \right\}$$
	is $\prod_{i\in I}\left\{\left(x_i,y_i\right)\in X_i\times X_i\mid x_i\leq_i y_i \right\}$, which is closed. 
\endproof 

\remark\label{r:powers of [0,1] are in PosComp}
	For every set $A$, $[0,1]^A$ (with the product order and product topology) is a compact ordered space.
\endremark 

If $X$ a Hausdorff space, then the diagonal of $X\times X$ is closed; as a consequence, we have the following.
\remark\label{r:hausdorff}
	Let $Y$ be a Hausdorff space, let $X$ be a topological space, and let $f,g\colon X\rightarrow Y$ be continuous functions. Then, $\{x\in X\mid f(x)=g(x)\}$ is closed.
\endremark

We can now prove the implication [\eqref{i:in-the-image}$\Rightarrow$\eqref{i:comp-pospace}] of Theorem~\ref{t:characterization fixed geometric}.
\theorem \label{t:X(A) is compact}

	For $A\in \MC$, $\Max(A)$ is a compact ordered space.
\endtheorem
\proof 
	$\Max(A)=\hom_{\MC}(A,[0,1])$ is a subset of $[0,1]^A$. 
	By Remark~\ref{r:powers of [0,1] are in PosComp}, $[0,1]^A$ (with the product order and product topology) is a compact ordered space. The topology on $\Max(A)$ coincides with the induced topology on $\Max(A)$ as a subspace of $[0,1]^A$; moreover, the order on $\Max(A)$ coincides with the order induced by $[0,1]^A$. 
	By Lemma~\ref{l:PosComp is hereditary}, it is enough to show that $\Max(A)$ is closed. The idea is that $\Max(A)$ is closed because it is defined by equations, which express the preservations of primitive operation symbols of $\MC$. 
	To make this precise, let $\mathcal{L}$ denote the set of primitive operation symbols of $\MC$. For each $h\in \mathcal{L}$, we denote with $\ar(h)$ the arity of $h$; moreover, we denote with $h_A$ the interpretation of $h$ in $A$, and by $h_{[0,1]}$ the interpretation of $h$ in $[0,1]$. For $a\in A$, we denote with $\pi_a\colon [0,1]^A\rightarrow [0,1]$ the projection onto the $a$-th coordinate (which is continuous). 
	We have
	$$
	\mld \Max(A)&=\hom_{\MC}(A,[0,1])=\\
	=\{x\colon A\rightarrow [0,1]\mid \forall h\in\mathcal{L}\ \forall a_1,\dots,a_{\ar(h)}\in A\\
	\ \ \ \ \ \ x(h_A(a_1,\dots,a_{\ar(h)}))= h_{[0,1]}(x(a_1), \dots, x(a_{\ar(h)}))\}=\\
	=\bigcap_{h\in \mathcal{L}}\bigcap_{a_1,\dots,a_{\ar(h)}\in A}\{x\colon A\rightarrow [0,1]\mid \\
	\ \ \ \ \ \ x(h_A(a_1,\dots,a_{\ar(h)}))= h_{[0,1]}(x(a_1), \dots, x(a_{\ar(h)}))\}=\\
	=\bigcap_{h\in \mathcal{L}}\bigcap_{a_1,\dots,a_{\ar(h)}\in A}\{x\in [0,1]^A\mid \pi_{h_A(a_1, \dots, a_{\ar(h)})}(x)= h_{[0,1]}(\pi_{a_1}(x), \dots, \pi_{a_{\ar(h)}}(x))\}.
	$$
	By Remark~\ref{r:all operations in V are mon and cont}, $h_{[0,1]}$ is continuous; therefore, the function from $[0,1]^A$ to $[0,1]$ which maps $x$ to $ h_{[0,1]}(\pi_{a_1}(x), \dots, \pi_{a_{\ar(h)}}(x))$ is continuous. Since $[0,1]$ is Hausdorff, by Remark~\ref{r:hausdorff}, $\{x\in [0,1]^A\mid \pi_{h_A(a_1, \dots, a_{\ar(h)})}(x)= h_{[0,1]}(\pi_{a_1}(x), \dots, \pi_{a_{\ar(h)}}(x))\}$ is closed.
\endproof

\subsection{The unit $\eta_X\colon X\rightarrow \Max\C(X)$ is injective}

We now turn to the proof of the implication [\eqref{i:comp-pospace}$\Rightarrow$\eqref{i:unit-iso}] in Theorem~\ref{t:characterization fixed geometric}, which states that, if $X$ is a compact ordered space, $\eta_X\colon X\rightarrow \Max\C(X)$ is an isomorphism. Our source of inspiration is \citep{HofmannLong}. The results we will obtain in the present section may be seen, essentially, as specific cases of the results available in \citep{HofmannLong}; nevertheless, for reasons of presentation, we provide independent proofs here.  In this subsection, we prove that, if $X$ is a compact ordered space, then $\eta_X\colon X\rightarrow \Max\C(X)$ is injective; this result is essentially due to L.\ Nachbin.

\definition
	For $X$ a partially ordered set, we call \emph{upper} an upward closed subset of $X$, and \emph{lower} a downward closed one.
	
\enddefinition
For $X$ a partial ordered set, and $x\in X$, we set $\downarrow x\coloneqq \{z\in X\mid z\leq x \}$ and $\uparrow x\coloneqq \{z\in X\mid x\leq z\}$. The following is well known.
\lemma\label{l:downarrow closed}
	Let $X$ be a compact ordered space, and let $x\in X$. Then $\downarrow x$ is the smallest closed lower subset of $X$ that contains $x$ and $\uparrow x$ is the smallest closed upper subset of $X$ that contains $x$.
\endlemma
\proof 
	Let us prove that $\downarrow x$ is closed. Set $D\coloneqq \{(u,v)\in X\times X\mid u\leq v \}$. $D$ is closed by definition of compact ordered space. Moreover, since any compact ordered space is Hausdorff, every point of $X$ is closed. Hence $D\cap (X\times \{x\})=\{(z,x)\mid z\in X: z\leq x \}$ is closed. Since the projection $\pi_1\colon X\times X\rightarrow X$ onto the first coordinate is closed, $\pi_1( \{(z,x)\mid z\in X: z\leq x \})=\{z\in X\mid z\leq x \}=\downarrow x$ is closed. Analogously for $\uparrow x$. The rest of the statement is straightforward to prove.
\endproof 

\proposition[Ordered version of Urysohn's Lemma]\label{p:Urysohn}
	Let $X$ be a compact ordered space, let $A$ be a closed lower subset, and let $B$ be a closed upper subset, with $A\cap B=\emptyset$. Then there exists a monotone and continuous function $\psi\colon X\rightarrow [0,1]$ such that, for every $x\in A$, $\psi(x)=0$, and, for every $x\in B$, $\psi(x)=1$.
\endproposition
\proof 
	See \citep[Chapter I, Theorem 1, p.\ 30]{Nachbin}.
\endproof 

\corollary\label{c:cor of urysohn}
	Let $X$ be a compact ordered space, and let $x,y\in X$ such that $x\not\geq y$. Then there exists a monotone and continuous function $\psi\colon X\rightarrow [0,1]$ such that $\psi(x)=0$ and $\psi(y)=1$.
\endcorollary
\proof 
	Set $\downarrow x\cap\uparrow y=\emptyset$. By Lemma~\ref{l:downarrow closed}, $\downarrow x$ is a closed lower subset and $\uparrow y$ is a closed upper subset. Therefore we may apply Proposition~\ref{p:Urysohn} with $A=\downarrow x$ and $B=\uparrow y$.
\endproof 
\corollary\label{c:cor of cor of urysohn}
	Let $X$ be a compact ordered space, and let $x,y\in X$. Suppose that, for every $\psi\colon X\rightarrow [0,1]$ monotone and continuous, $\psi(x)\leq \psi(y)$. Then $x\leq y$.
\endcorollary
A consequence of Corollary~\ref{c:cor of cor of urysohn} is the fact that every compact ordered space embeds into a power of $[0,1]$.

\proposition\label{p:eta injective}
	For every $X$ compact ordered space, $\eta_X$ is injective.
\endproposition
\proof 
	Let $x,y\in X$. Suppose $x\neq y$. Then, either $x\not\geq y$ or $y\not\geq x$. Suppose, without loss of generality, $x\not\geq y$. Then, by Corollary~\ref{c:cor of urysohn}, there exists $\psi\in \C(X)$ such that $\psi(x)=0$ and $\psi(y)=1$. Therefore, $(\eta_X(x))(\psi)=\ev_x(\psi)=\psi(x)=0\neq 1=\psi(y)=\ev_y(\psi)=(\eta_X(y))(\psi)$. Thus, $\eta_X(x)\neq \eta_X(y)$.
\endproof 

\subsection{The unit $\eta_X\colon X\rightarrow \Max\C(X)$ is surjective}
We continue the path that allows us to prove the implication [\eqref{i:comp-pospace}$\Rightarrow$\eqref{i:unit-iso}] in Theorem~\ref{t:characterization fixed geometric}: if $X$ is a compact ordered space, $\eta_X\colon X\rightarrow \Max\C(X)$ is an isomorphism. In this subsection, we prove that, if $X$ is a compact ordered space, then $\eta_X\colon X\rightarrow \Max\C(X)$ is surjective.

Let $X$ be a compact ordered space, and let $\Phi\colon \C (X)\rightarrow [0,1]$ be an MC-morphism, i.e.\ $\Phi\in \Max\C(X)$. The goal is to find $x\in X$ such that $\Phi=\ev_x$. For every $\psi\in \C(X)$, set $\Zm(\psi)\coloneqq \{x\in X\mid \psi(x)=0\}$. Moreover, set $\Zt(\Phi)\coloneqq \bigcap_{\psi\in \C(X):\Phi(\psi)=0}\Zm(\psi)$. We shall prove that $\Zt(\Phi)$ has a maximum element $x$, and that $\Phi=\ev_x$. Set $\A(\Phi)\coloneqq\bigcap_{\psi\in \C(X)}\{y\in X\mid \psi(y)\leq \Phi(\psi) \}$.

\lemma\label{l:A=phi}
	$\A(\Phi)= \Zt(\Phi)$.
\endlemma
\proof 
	Let us prove $(\subseteq)$. Let $y\in \A(\Phi)$. Then for every $\psi\in \C(X)$, we have $\psi(y)\leq \Phi(\psi)$. Therefore, if $\Phi(\psi)=0$, then $\psi(y)=0$, i.e., $\psi\in\Zm(\Phi)$. Hence, $y\in \Zt(\Phi)$.
	
	Let us prove $(\supseteq)$.  Let $x\in \Zt(\Phi)$. Let $\psi\in \C(X)$. We shall prove $\psi(x)\leq \Phi(\psi)$. Set $\psi'\coloneqq\psi\ominus \Phi(\psi)$. Then $\Phi(\psi')=\Phi(\psi)\ominus \Phi(\psi)=0$. Since $x\in \Zt(\Phi)$, we have $\psi'(x)=0$, i.e. $\psi(x)\ominus \Phi(\psi).$, i.e. $\psi(x)\leq \Phi(\psi)$.
\endproof 

\lemma\label{l:inequality A}
	\begin{enumerate}
		\item $\Zt(\Phi)$ is a closed lower subset of $X$.
		\item For every $\psi\in \C(X)$, $\sup_{y\in \A(\Phi)}\psi(y)\leq \Phi(\psi)$. 
		
	\end{enumerate}
	
\endlemma
\proof 
	\begin{enumerate}
		\item For every $\psi\in \C (X)$, $\Zm(\psi)$ is closed, hence $\Zt(\Phi)$ is closed. 
		
		Suppose $y\in \Zt(\Phi)$ and $x\leq y$. Then, for every $\psi\in \C(X)$ such that $\Phi(\psi)=0$, we have $y\in \Zm(\psi)$, i.e. $\psi(y)=0$. Since $\psi$ is monotone, we have $\psi(x)=0$, i.e. $x\in \Zm(\psi)$. Therefore, $x\in \Zt(\Phi)$.
		
		\item It follows from the fact that, for every $\psi\in \C(X)$ and $y\in \A(\Phi)$, $\psi(y)\leq \Phi(\psi)$.
	\end{enumerate}
\endproof 

The following is inspired by \citep[Proposition 6.12]{HofmannLong}.
\proposition
	Let $X$ be a compact ordered space and let $\Phi\colon \C(X)\rightarrow [0,1]$ be an MC-morphism. Then, for all $\psi\in \C(X)$, we have
	$$\Phi(\psi)= \sup_{y\in \Zt(\Phi)}\psi(y).$$
\endproposition
\proof 
	Let $\psi\in \C(X)$. We already know 
	$$\sup_{y\in \Zt(\Phi)}\psi(y)\stackrel{\text{Lem.~\ref{l:A=phi}}}{=}\sup_{y\in \A(\Phi)}\psi(y)\stackrel{\text{Lem.~\ref{l:inequality A}}}{\leq} \Phi(\psi).$$ Let us set $\lambda\coloneqq \sup_{y\in \Zt(\Phi)}\psi(y)$. We shall prove $\Phi(\psi)\leq \lambda$.  Let $\varepsilon>0$. We shall prove $\Phi(\psi)\leq \lambda+\varepsilon$. Set
	$$U\coloneqq\left\{y\in X\mid \psi(y)<\lambda+\varepsilon \right\}.$$
	Clearly, $U$ is open and, by definition of $\lambda$, $\Zt(\Phi)\subseteq U$. 
	Let $x\in X\setminus \Zt(\Phi)$. 
	There is some $\widetilde{\psi}\in \C(X)$ with $\Phi(\widetilde{\psi})=0$ and $\widetilde{\psi}(x)\neq 0$.
	Let $n\in \N$
	be such that $n(\widetilde{\psi}(x))\geq 1$.
	Set $\psi'\coloneqq\underbrace{\widetilde{\psi}\oplus\dots \oplus\widetilde{\psi}}_{n\text{ times}}$. 
	Then $\Phi(\psi')=0$ and $\psi'(x)=1$. 
	For every $\psi'\in \C(X)$ we set
	$$s(\psi')\coloneqq \{x\in X\mid \psi'(x)>1-\varepsilon \}.$$
	By the considerations above,
	$$X=U\cup\bigcup_{\psi' \in \C(X):\Phi(\psi')=0}s(\psi');$$
	since $X$ is compact, we find $\psi_1,\dots,\psi_n\in \C(X)$ with $\Phi(\psi_i)=0$ and
	$$X=U\cup s(\psi_1)\cup\dots\cup s(\psi_n).$$
	Therefore, for all $x\in X$, either $x\in U$, i.e., $\psi(x)<\lambda+\varepsilon$, or there exists $j\in \{1,\dots,n \}$ such that $x\in s(\psi_j)$, i.e., $\psi_j(x)>1-\varepsilon$.
	Hence,
	$$\psi\odot (1-\varepsilon)\leq \left(\left(\lambda\oplus\varepsilon\right)\odot (1-\varepsilon)\right)\lor \psi_1\lor \dots\lor \psi_n\leq \lambda\lor \psi_1\lor \dots\lor \psi_n.$$
	Therefore,
	$$\Phi(\psi)\odot (1-\varepsilon)\leq \lambda\lor \Phi(\psi_1)\lor\dots\lor \Phi(\psi_n)=\lambda\lor 0\lor \dots\lor 0=\lambda.$$
	Hence $\Phi(\psi)\leq(\Phi(\psi)\odot (1-\varepsilon))\oplus \varepsilon\leq \lambda\oplus \varepsilon$.
\endproof 
\definition
	A nonempty subset $A$ of a topological space $X$ is
	irreducible if $A\subseteq B\cup C$ for closed subsets $B$ and $C$ implies $A\subseteq B$ or $A \subseteq C$.
	
\enddefinition

\notation
	Let $X$ be a compact ordered space, and let $\tau$ be the topology on $X$. We denote with $\tau^\sharp$ the topology of all upper open subsets of $X$.
\endnotation
A subset $B\subseteq X$ is upper if, and only if, its complement $X\setminus B$ is lower; hence the closed subsets of $(X,\tau^\sharp)$ are precisely the closed lower subsets of $(X,\tau)$. 
\definition
	Let $X$ be a compact ordered space and let $\tau$ be the topology on $X$. We say that $A\subseteq X$ is \emph{$\sharp$-irreducible} if it is irreducible in $(X,\tau^\sharp)$, i.e., for all closed lower subsets $A_1,A_2$ of $(X,\tau)$ with $A\subseteq A_1\cup A_2$, one has $A\subseteq A_1$ or $A\subseteq A_2$.
	
\enddefinition

\proposition
	Let $X$ be a compact ordered space and $\Phi\colon \C(X)\rightarrow [0,1]$ an MC-morphism. Then
	\begin{enumerate}
		\item $\Zt(\Phi)\neq \emptyset$.
		\item $\Zt(\Phi)$ is $\sharp$-irreducible.
	\end{enumerate}
\endproposition
\proof 
	\begin{enumerate}
		\item $1=\Phi(1)=\sup_{x\in \Zt(\Phi)}1$.
		\item Let $A_1,A_2\subseteq X$ be closed upper subsets with $\Zt(\Phi)\subseteq A_1\cup A_2$. We shall prove that either $\Zt(\Phi)\subseteq A_1$ or $\Zt(\Phi)\subseteq A_2$. Suppose, by way of contradiction,  $x\in \Zt(\Phi)\setminus A_1$ and $y\in \Zt(\Phi)\setminus A_2$. Then $(\uparrow x )\cap A_1=\emptyset$ and $(\uparrow y)\cap A_2=\emptyset$. By Proposition~\ref{p:Urysohn}, there exist $\psi_1,\psi_2\in \C(X)$ such that, 
		\begin{enumerate}
			\item $\psi_1(x)=1$ and, for all $z\in A_1$, $\psi_1(z)=0$.
			\item $\psi_2(y)=1$ and, for all $z\in A_2$, $\psi_2(z)=0$.
		\end{enumerate}
		Then, for all $z\in \Zt(\Phi)$, $(\psi_1\odot \psi_2)(z)=\psi_1(z)\odot \psi_2(z)=0$, since either $\psi_1(z)=0$ or $\psi_2(z)=0$.
		
		Hence, $\Phi(\psi_1\odot \psi_2)=\sup_{z\in \Zt(\Phi)}\psi_1(z)\odot \psi_2(z)=0$.
		
		But, also, 
		$$
		\mld \Phi(\psi_1\odot \psi_2)&=\Phi(\psi_1)\odot \Phi(\psi_2)=\\
		=\left(\sup_{z\in \Zt(\Phi)}\psi_1(z)\right)\odot\left(\sup_{z\in \Zt(\Phi)}\psi_2(z)\right)\geq\\
		\geq \psi_1(x)\odot \psi_2(y)=\\
		=1\odot 1=\\
		=1.
		$$
		This is a contradiction.
	\end{enumerate}
\endproof 

\definition
	We say that a topological space $(X,\tau)$ is \emph{sober} if, for every irreducible closed set $C$, there exists a unique $x\in X$ such that the closure of $\{x\}$ is $C$.
	
\enddefinition
\proposition\label{p:sharp is sober}
	Let $(X,\leq, \tau)$ be a compact ordered space. Then $(X,\tau^\sharp)$ is sober.
\endproposition
\proof 
	See Proposition VI.6.11 in \citep{ContLattDom}.
\endproof 

\remark\label{r:closure in sharp}
	Let $(X,\leq,\tau)$ be a compact ordered space, and let $x\in X$. The closure of $\{x\}$ in $(X,\tau^\sharp)$ is $\downarrow x$.
\endremark 

\proposition\label{p:irreducible is principal}
	Let $X$ be a compact ordered space, and let $A\subseteq X$ be an irreducible closed lower subset of $X$. Then there exists a unique $x\in X$ such that $A=\downarrow x$.
\endproposition
\proof 
	By Proposition~\ref{p:sharp is sober} and Remark~\ref{r:closure in sharp}.
\endproof 

\theorem 
	Let $X$ be a compact ordered space, and let $\Phi\colon \C(X)\rightarrow [0,1]$ be a map that preserves the operations $\lor,\land,\oplus,\odot$ and every constant $\lambda\in [0,1]$. Then, there exists a unique $x\in X$ such that, for every $\psi\in \C(X)$, $\Phi(\psi)=\psi(x)$, i.e. $\Phi=\ev_x$.
\endtheorem
\proof 
	For every $x\in X$, $\eta_X(x)=\ev_x$. Hence, uniqueness follows from injectivity of $\eta_X$, which was established in Proposition~\ref{p:eta injective} and which---we recall---was a consequence of Corollary~\ref{c:cor of urysohn}.
	Concerning existence, by Proposition~\ref{p:irreducible is principal}, there exists $x\in X$ such that $\Zt(\Phi)=\downarrow x$. Then
	$$\Phi(\psi)= \sup_{z\in \Zt(\Phi)}\psi(z)=\psi(x).$$
\endproof 
\corollary\label{c:eta surjective}
	If $X$ is a compact ordered space, the map $\eta_X\colon X\rightarrow \Max\C(X)$ is surjective.
\endcorollary
We may now conclude the proof of Theorem~\ref{t:characterization fixed geometric}, which asserted, for $X$ a preordered topological space, the equivalence of the following conditions.
\begin{enumerate}
	\item The unit $\eta_X\colon X\rightarrow \Max\C(A)$ is an isomorphism.
	\item There exists an MC-algebra $A$ such that $X$ and $\Max(A)$ are isomorphic preordered topological spaces.
	\item $X$ is a compact ordered space.
\end{enumerate}

\proof [of Theorem~\ref{t:characterization fixed geometric}]
	\noindent	[\eqref{i:unit-iso}$\Rightarrow$\eqref{i:in-the-image}] By Remark~\ref{r:specialization}.
	
	\noindent [\eqref{i:in-the-image}$\Rightarrow$\eqref{i:comp-pospace}] By Theorem~\ref{t:X(A) is compact}.
	
	\noindent [\eqref{i:comp-pospace}$\Rightarrow$\eqref{i:unit-iso}]  By Proposition~\ref{p:eta injective}, $\eta_X$ is injective. By Corollary~\ref{c:eta surjective}, $\eta_X$ is surjective. 	Every continuous map between compact Hausdorff spaces is closed, and every closed bijective continuous map between topological spaces is a homeomorphism. We are left to show that $\eta_X$ reflects the order, i.e., for every $x,y\in X$, if $\eta_X(x)\leq \eta_X(y)$, then $x\leq y$. If $\eta_X(x)\leq \eta_X(y)$, then, for every $\psi\in \C (X)$, $(\eta_X(x))(\psi)\leq (\eta_X(y))(\psi)$, i.e., $\psi(x)\leq \psi(y)$. By Corollary~\ref{c:cor of cor of urysohn}, $x\leq y$.
\endproof

\section{Fixed objects on the algebraic side: the goal}
\definition
	Let $A\in\MC$ and $x,y\in A$. We set 
	$$\uparrow_x^y\coloneqq \{\lambda\in[0,1]\mid y\leq x\oplus \lambda\};$$
	$$\du(x,y)\coloneqq \inf\uparrow_x^y;$$
	and
	$$\dist(x,y)\coloneqq \max\{\du(x,y), \du(y,x) \}.$$
	
\enddefinition
On $[0,1]$, $\dist^{\uparrow}(x,y)=(y-x)^+$ (where $z^+\coloneqq\max\{z,0\}$), and $\dist(x,y)=\lvert y-x\rvert$. If $X$ is a set, and $L$ is an MC-subalgebra of $[0,1]^X$, then, on $L$, $\du(f,g)=\sup_{x\in X}(g(x)-f(x))^+$, and $\dist$ coincides with the sup metric. We mention that, in $[0,1]$, $(y-x)^+$ coincides with $y\ominus x$.

\definition
	Let $A\in\MC $. We say that $A$ is \emph{archimedean} if, for all $x,y\in A$, $$\dist(x,y)=0 \Rightarrow x=y.$$
	
\enddefinition
The idea---as we will see---is that $A\in \MC$ is archimedean if, and only if, $A$  is an MC-subalgebra of $[0,1]^X$, for some set $X$. For now, we have the following.
\remark\label{r:subalgebra then arch}
	If $X$ is a set, and $L$ is an MC-subalgebra of $[0,1]^X$, then $L$ is archimedean. Indeed,  $\dist$ coincides with the sup metric, that satisfies the implication $\dist(x,y)=0 \Rightarrow x=y$.
\endremark 

In Definition~\ref{d:Cauchy} below, we define Cauchy sequences, convergence, and Cauchy completeness. These definitions are standard; anyway, one should pay attention to the fact that $\dist$ is not required to be a metric, because $\dist(x,y)=0\Rightarrow x=y$ might fail.
\definition\label{d:Cauchy}
	Let $A\in \MC$, let $(a_n)_{n\in \N}$ be a sequence in $A$, and let $a\in A$. We say that $(a_n)_{n\in \N}$ is a \emph{Cauchy sequence} if, for all $\varepsilon>0$, there exists $k\in \N$ such that, for all $n,m\geq k$, $\dist(a_{n},a_{m})<\varepsilon$. We say that  $(a_n)_{n\in \N}$ \emph{converges to $a$}, or that $a$ is a limit of $(a_n)_{n\in \N}$, if, for every $\varepsilon>0$, there exists $n\in \N$ such that, for all $m\geq n$, $\dist(a_m,a)<\varepsilon$. We say that $(a_n)_{n\in \N}$ \emph{converges} if there exists $b\in A$ such that $(a_n)_{n\in \N}$ converges to $b$. We say that $A$ is \emph{Cauchy complete} if every Cauchy sequence in $A$ converges.
	
\enddefinition
We remark that a sequence may have more than one limit, since $\dist$ is not a metric.

\remark
	On $[0,1]$ the concepts of Cauchy sequence and convergence to an element in Definition~\ref{d:Cauchy} coincide with the usual ones with respect to the euclidean distance. In particular, $[0,1]$ is Cauchy complete. If $X$ is a set, and $L$ is an MC-subalgebra of $[0,1]^X$, the concepts of Cauchy sequence and convergence to an element in Definition~\ref{d:Cauchy} coincide with the usual ones with respect to the sup metric.
\endremark 

Our next goal is to prove the following.
\theorem \label{t:characterization algebraic fixed objects}
	Let $A\in \MC$. The following conditions are equivalent.
	\begin{enumerate}
		\item\label{i:unit-iso-alg} The unit $\varepsilon_A\colon A\rightarrow \C\Max(A)$ is an isomorphism.
		\item\label{i:in-the-image-alg} There exists a preordered topological space $X$ such that $A$ and $\C(X)$ are isomorphic MC-algebras.
		\item\label{i:arch-cc} $A$ is archimedean and Cauchy complete.
	\end{enumerate}
\endtheorem

The current section and the following two are intended to prove Theorem~\ref{t:characterization algebraic fixed objects} above.

\remark\label{r:easy}
	We give here the proofs of the implications [\eqref{i:unit-iso-alg}$\Rightarrow$\eqref{i:in-the-image-alg}] and [\eqref{i:in-the-image-alg}$\Rightarrow$\eqref{i:arch-cc}] of Theorem~\ref{t:characterization algebraic fixed objects}.
	
	\noindent [\eqref{i:unit-iso-alg}$\Rightarrow$\eqref{i:in-the-image-alg}] Take $X=\Max(A)$.
	
	\noindent [\eqref{i:in-the-image-alg}$\Rightarrow$\eqref{i:arch-cc}] $A$ is archimedean because $\C(X)$ is archimedean, by Remark~\ref{r:subalgebra then arch}. Let us prove that $\C(X)$ (and hence $A$) is Cauchy complete. Let $(f_n)_{n\in \N}$ be a sequence in $\C(X)$ which is Cauchy with respect to the metric $\dist$. Then, there exists a function $f\colon X\rightarrow [0,1]$ such that $f_n$ converges to $f$ uniformly. It is well known that the uniform limit of a sequence of continuous functions is continuous. Since, for all $n\in \N$, $f_n$ is continuous, $f$ is continuous. Let us prove that $f$ is monotone. Let $x,y\in X$ with $x\leq y$. For all $n\in \N$, $f_n$ is monotone. Therefore, $f(x)=\lim_{n\rightarrow \infty}f_n(x)\leq \lim_{n\rightarrow \infty}f_n(y)=f(y)$.
\endremark

We are left to prove the implication [\eqref{i:arch-cc}$\Rightarrow$\eqref{i:in-the-image-alg}] of Theorem~\ref{t:characterization algebraic fixed objects}, i.e., that, for every archimedean Cauchy complete MC-algebra $A$, the unit $\varepsilon_A\colon A\rightarrow \C\Max(A)$ is an isomorphism. This implication is a consequence of the following two theorems, whose proofs we conclude, respectively,  at the ends of the following two sections.

\theorem \label{t:arch iff}
	Let $A\in \MC$. The following conditions are equivalent.
	\begin{enumerate}
		\item\label{i:arch} $A$ is archimedean.
		\item\label{i:enough-morph} For every $x,y\in A$ with $x\neq y$, there exists an MC-morphism $\varphi\colon A\rightarrow [0,1]$ such that $\varphi(x)\neq \varphi(y)$.
		\item\label{i:subalg} There exists a set $X$ such that $A$ is an MC-subalgebra of $[0,1]^X$.
		\item\label{i:unit-inj} The unit $\varepsilon_A\colon A\rightarrow \C\Max(A)$ is injective.
	\end{enumerate}
\endtheorem
\theorem \label{t:Cauchy complete iff}
	Let $A\in\MC$. The following conditions are equivalent.
	\begin{enumerate}
		\item $A$ is Cauchy complete.
		\item The unit $\varepsilon_A\colon A\rightarrow \C\Max(A)$ is surjective.
	\end{enumerate}
\endtheorem

\section{$A$ is archimedean if, and only if, the unit $\varepsilon_A$ is injective}
The aim of this section is to prove Theorem~\ref{t:arch iff} above. Some of the implications between the four conditions in Theorem~\ref{t:arch iff} are relatively easy to prove, and we collect their proofs in the following remark.

\remark[Part of the proof of Theorem~\ref{t:arch iff}]\label{r:proof arch}
	We prove some of the implications in Theorem~\ref{t:arch iff}.
	
	\noindent[\eqref{i:unit-inj}$\Rightarrow$\eqref{i:subalg}] Immediate.
	
	\noindent [\eqref{i:subalg}$\Rightarrow$\eqref{i:enough-morph}] Let $f,g\in A\subseteq [0,1]^X$ be such that $f\neq g$. Then, there exist $z\in X$ such that $f(z)\neq g(z)$. Set
	\begin{align*}
	\varphi\colon A&\longrightarrow [0,1]\\
	h&\longmapsto h(z).
	\end{align*} The function $\varphi\colon A\rightarrow [0,1]$ is clearly an MC-morphism. $\varphi(f)=f(z)\neq g(z)=\varphi(g)$.
	\noindent[\eqref{i:enough-morph}$\Rightarrow$\eqref{i:unit-inj}] Let $x,y\in A$ be such that $x\neq y$. We shall prove $\varepsilon_A(x)\neq \varepsilon_A(y)$. By hypothesis, there exists an MC-morphism $\varphi\colon A\rightarrow [0,1]$ such that $\varphi(x)\neq \varphi(y)$. Note that $\varphi\in \Max(A)$. We have $(\varepsilon_A(x))(\varphi)=\ev_x(\varphi)=\varphi(x)\neq \varphi(y)=\ev_y(\varphi)=(\varepsilon_A(y))(\varphi)$. This proves $\varepsilon_A(x)\neq \varepsilon_A(y)$.
	
	\noindent  [\eqref{i:subalg}$\Rightarrow$\eqref{i:arch}] By Remark~\ref{r:subalgebra then arch}.
\endremark 
In Remark~\ref{r:proof arch} we have proved that \eqref{i:enough-morph}, \eqref{i:subalg} and \eqref{i:unit-inj} in Theorem~\ref{t:arch iff} are equivalent, and that any of these conditions implies \eqref{i:arch}. What is missing to prove Theorem~\ref{t:arch iff} is the implication [\eqref{i:arch}$\Rightarrow$\eqref{i:enough-morph}] (for example), i.e., that if $A$ is archimedean, then, for every $x\neq y\in A$, there exists an MC-morphism $\varphi\colon A\rightarrow [0,1]$ such that $\varphi(x)\neq \varphi(y)$. The proof of this last implication will take us some effort. The idea is that, since $\MC$ is a \textit{finitary} variety, we may apply the Subdirect Representation Theorem. The Subdirect Representation Theorem ensures, for every $A\in \MC$, the existence of an injective MC-morphism $\iota\colon A\hookrightarrow \prod_{i\in I}A_i$, where, for all $i\in I$, $A_i$ is subdirectly irreducible. We will show that every subdirectly irreducible algebra $B$ consists essentially of the set $[0,1]$ together with some additional elements, each of which lies ``just above'' or ``just below'' one particular $\lambda\in [0,1]$. This allows us to define the morphism ${\rm ess}\colon B\rightarrow [0,1]$ that ``kills the infinitesimals''. We will show that the morphism ${\rm ess}$ preserves $\dist$. Composing the map $\iota\colon A\hookrightarrow \prod_{i\in I}A_i$ with the morphisms ${\rm ess}_i\colon A_i\rightarrow [0,1]$, we obtain a map $h\colon A\rightarrow\prod_{i\in I}[0,1]$ that correctly translates the function $\dist$ on $A$ to the sup metric on $\prod_{i\in I}[0,1]$. Then, if, $x\neq y\in A$, and $A$ is archimedean, we have $\dist(x,y)\neq0$, and therefore the sup distance between $h(x)$ and $h(y)$ is not zero. From this we conclude that there is a morphism $\varphi\colon A\rightarrow [0,1]$ such that $\varphi (x)\neq \varphi(y)$.

\subsection{Subdirectly irreducible MC-algebras}
We start the path that will take us to the proof that every archimedean MC-algebra admits enough morphisms towards $[0,1]$, i.e., the implication [\eqref{i:arch}$\Rightarrow$\eqref{i:enough-morph}] in Theorem~\ref{t:arch iff}.
The main goal of this subsection is to prove the following.
\theorem \label{t:irreducible is so}
	Let $A\in\MC$ be subdirectly irreducible. Then, for all $x\in A$, and $\lambda\in [0,1]$, $x\leq \lambda$ or $\lambda\leq x$.
\endtheorem

To a reader that has familiarity with abelian lattice-ordered groups or  MV-algebras, Theorem~\ref{t:irreducible is so} might sound analogous  to the well-known result that every subdirectly irreducible abelian lattice-ordered group (or MV-algebra) is totally ordered. Note, however, that Theorem~\ref{t:irreducible is so} does not say ``every subdirectly irreducible $A\in \MC$ is totally ordered''; in fact, it says something something weaker. The author does not know whether every subdirectly irreducible $A\in \MC$ is totally ordered; nevertheless, for our purposes, Theorem~\ref{t:irreducible is so} is enough.

The idea for the proof of Theorem~\ref{t:irreducible is so} is the following: if, by contraposition, $x\not\leq \lambda$ and $\lambda\not\leq x$, then we could construct two nonminimal congruences whose intersection is minimal, which shows that $A$ is not subdirectly irreducible. The idea is that these two congruences are the congruences generated, respectively, by $x\sim x\lor \lambda$ and $x\sim x\land \lambda$.

Our convention is that we do not consider the trivial algebra (i.e, with exactly an element) as subdirectly irreducible.

\notation
	Let $A\in \MC$, and let $\mathcal{L}$ denote the language of $\MC$. For $A\in \MC$, let $A^\partial$ denote the $\mathcal{L}$-algebra that shares the same underlying set with $A$, and which is such that $\lor_{A^\partial}=\land_A$, $\land_{A^\partial}=\lor_A$, $\oplus_{A^\partial}=\odot_A$, $\odot_{A^\partial}=\oplus_A$, and, for every $\lambda\in [0,1]$, $\lambda_{A^\partial}=(1-\lambda)_A$. We call $A^\partial$ the \emph{order-dual algebra} of $A$.
\endnotation
Roughly speaking, the dual operation $\gamma^\partial$ of an operation $\gamma$ of arity $I$  is given by $\gamma^{\partial}((x_i)_{i\in I})\coloneqq1-\gamma((1-x_i)_{i\in I}) $ (this makes sense in $[0,1]$).

We will use the concept of order-dual algebra only to shorten some proofs. This is made possible by the following lemma.

\lemma\label{l:dual algebra}
	Let $A\in \MC$. Then, $A^\partial\in \MC$.
\endlemma
\proof 
	The only nontrivial part is showing that Axioms \eqref{i:lor}, \eqref{i:land}, \eqref{i:oplus}, \eqref{i:odot} ``dualize''. But this holds because, for every $a,b,c\in [0,1]$, we have the following properties.
	\begin{enumerate}
		\item $a\land b=c\Leftrightarrow (1-a)\lor (1-b)=1-c$.
		\item $a\lor b=c\Leftrightarrow (1-a)\land (1-b)=1-c$.
		\item $a\oplus b=c\Leftrightarrow (1-a)\odot (1-b)=1-c$.
		\item $a\odot b=c\Leftrightarrow (1-a)\oplus (1-b)=1-c$.
	\end{enumerate} 
\endproof

\remark
	$\left(A^\partial\right)^\partial=A$.
\endremark 

\lemma
	Let $A\in \MC$. The following properties hold for all $a,b,c,a',b'\in A$.
	\begin{enumerate}
		\item\label{i:transl-inv-oplus} $a\leq b\Rightarrow a\oplus c\leq b\oplus c$.
		\item\label{i:transl-inv-odot} $a\leq b\Rightarrow a\odot c\leq b\odot c$.
		\item\label{i:oplus-adds} $a\leq a\oplus b$.
		\item\label{i:odot-subtracts} $a\geq a\odot b$.
		\item\label{i:oplus-monotone} If $a\leq a'$ and $b\leq b'$, then $a\oplus b\leq a'\oplus b'$.
		\item\label{i:odot-monotone} If $a\leq a'$ and $b\leq b'$, then $a\odot b\leq a'\odot b'$.
	\end{enumerate}
\endlemma
\proof 
	\noindent \eqref{i:transl-inv-oplus} Suppose $a\leq b$. Then $a\oplus c=(a\land b)\oplus c=(a\oplus c)\land (b\oplus c)$.\\
	\noindent \eqref{i:oplus-adds}  From $0\leq b$ we obtain $a\oplus 0\leq a\oplus b$. Since $a\oplus 0=a$, $a\leq a\oplus b$.\\
	\noindent \eqref{i:oplus-monotone} This follows from \eqref{i:transl-inv-oplus}, taking into account that $\oplus$ is commutative.\\
	\noindent Items \eqref{i:transl-inv-odot}, \eqref{i:odot-subtracts} and \eqref{i:odot-monotone} in $A$ coincide respectively with items \eqref{i:transl-inv-oplus}, \eqref{i:oplus-adds} and \eqref{i:oplus-monotone} in the order-dual algebra of $A$.
\endproof

\lemma\label{l:sim+}
	Let $A\in \MC$ and $x\in A$. For $a,a'\in A$, set $ a\sim_\oplus^x a'$ if, and only if, there exist $n,m\in \N$ such that 
	$$ a\oplus(\underbrace{x\oplus\cdots\oplus x}_{n\text{ times}})\geq a'$$
	$$a'\oplus(\underbrace{x\oplus\cdots\oplus x}_{m\text{ times}})\geq a.$$ 	
	Then $\sim_\oplus^x$ is a congruence.
\endlemma
\proof 
	We first prove that $\sim_\oplus^x$ is an equivalence relation. The relation $\sim_\oplus^x$ is trivially reflexive and symmetric. To prove transitivity, suppose $a\sim_\oplus^x b\sim_\oplus^x c$. Then 
	$$(a\oplus(\underbrace{x\oplus\cdots\oplus x}_{n\text{ times}}))\oplus(\underbrace{x\oplus\cdots\oplus x}_{n'\text{ times}})\geq b\oplus (\underbrace{x\oplus\cdots\oplus x}_{n'\text{ times}})\geq c.$$
	Analogously for the other inequality.
	
	Suppose $a\sim_\oplus^x a'$ and $b\sim_\oplus^x b'$.
	
	\begin{enumerate}
		\item We shall prove $a\lor b\sim_\oplus^x a'\lor b'$.			
		$$(a\lor b)\oplus(\underbrace{x\oplus\cdots\oplus x}_{\max\{n,n'\}\text{ times}})\geq  a'\lor b'.$$			
		Analogously for the other inequality.			
		\item We shall prove $a\land b\sim_\oplus^x a'\land b'$.			
		$$(a\land b)\oplus(\underbrace{x\oplus\cdots\oplus x}_{\max\{n,n'\}\text{ times}})=(a\oplus(\underbrace{x\oplus\cdots\oplus x}_{\max\{n,n'\}\text{ times}}))\land (b\oplus(\underbrace{x\oplus\cdots\oplus x}_{\max\{n,n'\}\text{ times}})) \geq a'\land b'.$$
		Analogously for the other inequality.
		\item We shall prove $a\oplus b\sim_\oplus^x a'\oplus b'$.			
		$$(a\oplus b)\oplus(\underbrace{x\oplus\cdots\oplus x}_{(n+n')\text{ times}})=(a\oplus(\underbrace{x\oplus\cdots\oplus x}_{n\text{ times}}))\oplus (b\oplus(\underbrace{x\oplus\cdots\oplus x}_{n'\text{ times}})) \geq a'\oplus b'.$$
		Analogously for the other inequality.
		
		\item We shall prove $a\odot b\sim_\oplus^x a'\odot b'$.			 
		$$(a\odot b)\oplus(\underbrace{x\oplus\cdots\oplus x}_{(n+n')\text{ times}})\geq(a\oplus(\underbrace{x\oplus\cdots\oplus x}_{n\text{ times}}))\odot (b\oplus(\underbrace{x\oplus\cdots\oplus x}_{n'\text{ times}})) \geq a'\odot b'.$$
		Analogously for the other inequality.
	\end{enumerate}
\endproof 

\lemma\label{l:sim-}
	Let $A\in \MC$, $x\in A$. 
	
	For $b,b'\in A$, set $ b\sim_\odot^x b'$ if, and only if there exists $n,m\in \N$ such that 
	$$ b\odot(\underbrace{x\odot\cdots\odot x}_{n\text{ times}})\leq b',$$
	$$b'\odot(\underbrace{x\odot\cdots\odot x}_{m\text{ times}})\leq b.$$ 
	Then $\sim_\odot^x$ is a congruence.
\endlemma
\proof 
	The relation $\sim_\odot^x$ is a congruence in $A^\partial$, by Lemma~\ref{l:sim+}. This implies that $\sim_\odot^x$ is a congruence in $A$. 
\endproof 

We call \emph{minimal congruence} (or \emph{trivial congruence}) on an algebra $A$ the smallest congruence, which is $\Delta_A\coloneqq\{(a,a)\mid a\in A \}$.

We are now ready for the proof of Theorem~\ref{t:irreducible is so}.

\proof[of Theorem~\ref{t:irreducible is so}]
	Suppose, by way of contradiction, $x\not\leq \lambda$ and $\lambda\not\leq x$. Then, $x\ominus \lambda\neq 0$ and $x\oplus (1-\lambda)\neq 1$. Indeed, $x\leq (x\ominus \lambda)\oplus \lambda$, and $x\geq (x\oplus (1-\lambda))\ominus (1-\lambda)$.
	
	For $a,b\in A$, set $a\sim_\oplus b$ if, and only if, there exists $n,m\in \N$ such that 
	$$ a\oplus(\underbrace{(x\ominus\lambda)\oplus\cdots\oplus (x\ominus \lambda)}_{n\text{ times}})\geq b,$$
	$$b\oplus(\underbrace{(x\ominus \lambda)\oplus\cdots\oplus (x\ominus \lambda)}_{m\text{ times}})\geq a.$$ 
	By Lemma~\ref{l:sim+}, $\sim_\oplus$ is a congruence.
	
	For $a,b\in A$, set $a\sim_\odot b$ if, and only if, there exists $n,m\in \N$ such that 
	$$ a\odot(\underbrace{(x\oplus(1-\lambda))\odot\cdots\odot (x\oplus (1-\lambda))}_{n\text{ times}})\leq b,$$
	$$b\odot(\underbrace{(x\oplus(1-\lambda))\odot\cdots\odot (x\oplus(1-\lambda))}_{m\text{ times}})\leq a.$$ 
	By Lemma~\ref{l:sim-},	$\sim_\odot$ is a congruence.
	
	Since $x\ominus \lambda\neq 0$ and $x\ominus\lambda\sim_\oplus0$, $\sim_\oplus$ is not the minimal congruence. Since $x\oplus (1-\lambda)\neq 1$ and $x\oplus(1-\lambda)\sim_\oplus 1$, $\sim_\oplus$ is not the minimal congruence.
	
	We claim that the congruence $\sim_\oplus\cap\sim_\odot$ is the minimal one. Indeed, let us take $a,b\in A$ such that $a\sim_\oplus b$ and $a\sim_\odot b$. 
	
	Then $$a\leq b\oplus(\underbrace{(x\ominus \lambda)\oplus\cdots\oplus (x\ominus \lambda)}_{m\text{ times}}).$$ 
	and
	$$ a\odot(\underbrace{(x\oplus(1-\lambda))\odot\cdots\odot (x\oplus (1-\lambda))}_{n\text{ times}})\leq b.$$
	Then
	$$
	\mld  a&=a\land( b\oplus(\underbrace{(x\ominus \lambda)\oplus\cdots\oplus (x\ominus \lambda)}_{m\text{ times}}))\leq\\
	\leq (a\odot(\underbrace{(x\oplus(1-\lambda))\odot\cdots\odot (x\oplus (1-\lambda))}_{n\text{ times}}))\lor b=b.
	$$
	Analogously, $b\leq a$, and therefore $a=b$. Therefore $\sim_\oplus\cap\sim_\odot$ is the minimal congruence. This is a contradiction; indeed, the minimal congruence on a subdirectly irreducible algebra is $\land$-irreducible in the lattice of congruences, and therefore the intersection of nonminimal congruence is not minimal in a subdirectly irreducible algebra.
\endproof 

\subsection{The morphism ${\rm ess}$ from a subdirectly irreducible MC-algebra to $[0,1]$ that kills infinitesimals}
In Theorem~\ref{t:irreducible is so} we proved that, if $A\in\MC$ is subdirectly irreducible, then, for all $x\in A$, and $\lambda\in [0,1]$, $x\leq \lambda$ or $\lambda\leq x$. The intuition is that $A$ consists essentially of the set $[0,1]$, together with some infinitesimals, each of which lies ``just above'' or ``just below'' one particular $\lambda\in [0,1]$. In this subsection we show that one can define the function ${\rm ess}\colon A\rightarrow [0,1]$ that ``kills the infinitesimals'' and that this map is an MC-morphism (see Theorem~\ref{t:morphism to [0,1]}).

\lemma\label{l:trivial if and only if 0=1}
	Let $A\in\MC$. Let $\alpha,\beta\in[0,1]$ be such that $\alpha\neq \beta$. The following conditions are equivalent.
	\begin{enumerate}
		\item $A$ is trivial.
		\item $\alpha_A=\beta_A$.
	\end{enumerate}
\endlemma
\proof 
	If $A$ is trivial, then, clearly, $\alpha_A=\beta_A$. Suppose $\alpha_A=\beta_A$ and suppose, without loss of generality, $\alpha<\beta$. Then, in $A$, $\beta\ominus \alpha=\beta\ominus\beta=0$. Then, for every $n\in\N$, in $A$, $0=\underbrace{(\beta\ominus \alpha)\oplus \dots\oplus (\beta\ominus\alpha)}_{n \text{ times}}$. Let $n$ be big enough so that, in $[0,1]$, $\underbrace{(\beta\ominus \alpha)\oplus \dots\oplus (\beta\ominus\alpha)}_{n \text{ times}}=1$. Then, in $A$, $0_A=\underbrace{(\beta\ominus \alpha)\oplus \dots\oplus (\beta\ominus\alpha)}_{n \text{ times}}=1_A$. Since $0_A$ is the bottom and $1_A$ is the top of $A$, for every $x,y\in A$, $x=y$.
\endproof 
\definition
	Let $A\in \MC$ and $x\in A$. We set 
	$$I_x\coloneqq\{\lambda\in [0,1]\mid \lambda\leq x \};$$
	$$S_x\coloneqq\{\lambda\in [0,1]\mid x\leq\lambda \};$$
	$$\essinf x\coloneqq \sup I_x;$$
	$$\esssup x\coloneqq \inf S_x.$$
	
\enddefinition
\remark
	Let $A\in \MC$ and $x\in A$. Then, $\essinf x$, calculated in $A$, is $\esssup x$, calculated in $A^{\partial}$.
\endremark 

\remark\label{r:essinf leq esssup}
	$I_x$ is an initial segment of $[0,1]$, $0\in I_x$, $S_x$ is a final segment of $[0,1]$, and $1\in S_x$. In addition, we note the following facts.
	\begin{enumerate}
		\item If $A$ is trivial, then $I_x=[0,1]=S_x$, and thus $\essinf x=1$ and $\esssup x=0$.
		\item If $A$ is not trivial, then $I_x\cap S_x$ has at most one element, because otherwise there would exists $\alpha\neq\beta$ in $[0,1]$ such that, in $A$, $\alpha\leq x\leq \beta\leq x\leq \alpha$, and hence $\alpha_A=\beta_A$; but, since $A$ is non-trivial, this is not possible by Lemma~\ref{l:trivial if and only if 0=1}. Hence, if $A$ is not trivial,
		$$\essinf x\leq \esssup x.$$
	\end{enumerate} 
\endremark

\lemma\label{l:then essinf=esssup}
	Let $A\in\MC$ be nontrivial. If $x\in A$ is such that, for all $\lambda\in [0,1]$, $x\leq \lambda$ or $\lambda\leq x$, then
	$$\essinf x=\esssup x.$$
	In particular, if $A$ is subdirectly irreducible, then, for all $x\in X$, $\essinf x=\esssup x$.
\endlemma
\proof 
	By Remark~\ref{r:essinf leq esssup}, $\essinf x\leq \esssup x$. Since, by hypothesis, $I_x\cup S_x=[0,1]$, we conclude $ \essinf x=\esssup x$. If $A$ is subdirectly irreducible, then, by Theorem~\ref{t:irreducible is so}, for all $x\in A$ and $\lambda\in [0,1]$, $x\leq \lambda$ or $\lambda\leq x$.
\endproof 
\notation
	For $A\in \MC$ and $x\in A$, if $\essinf x=\esssup x$, we set
	$$\ess{x} \coloneqq \essinf x=\esssup x.$$ 
\endnotation

\theorem \label{t:morphism to [0,1]}
	Let $A\in\MC$ be such that, for all $x\in A$, $\essinf x=\esssup x$ (this holds, in particular, if $A$ is subdirectly irreducible). Then, the  function
	\begin{align*}
	{\rm ess}\colon \MC&\longrightarrow [0,1]\\
	x&\longmapsto\ess x
	\end{align*}
	is a surjective MC-morphism.
\endtheorem
\proof 
	For all $x\in A$, recall: $I_x\coloneqq\{\lambda\in [0,1]\mid \lambda\leq x \}$, $S_x\coloneqq\{\lambda\in [0,1]\mid x\leq\lambda \}$, and $\ess x=\sup I_x=\inf S_x$. For every constant symbol $\lambda\in [0,1]$, $\ess$ clearly preserves $\lambda$. Let $x,y\in A$ and let $\otimes$ denote any operation amongst $\{\lor, \land ,\oplus, \odot \}$. We shall show $\ess(x\otimes y)=\ess x\otimes\ess y$. For every $U,W\subseteq [0,1]$, set $U\otimes W\coloneqq\{\alpha\otimes\beta\mid \alpha\in U,\beta\in W \}$. Since $\otimes\colon [0,1]^2\rightarrow [0,1]$ is continuous, for every nonempty $U,V\subseteq [0,1]$, $\sup (U\otimes W)=(\sup U)\otimes (\sup W)$ and $\sup (U\otimes W)=(\sup U)\otimes (\sup W)$. We shall show $\ess(x\otimes y)=\ess x\otimes\ess y$. Now, 
	$$\ess(x\otimes y)=\sup I_{x\otimes y},$$
	$$\ess(x\otimes y)=\inf S_{x\otimes y},$$
	$$\ess x\otimes\ess y=(\sup I_x)\otimes (\sup I_y)=\sup(I_x\otimes I_y),$$ $$\ess x\otimes\ess y=(\inf S_x)\otimes (\inf S_y)=\inf(S_x\otimes S_y).$$
	Let us take $\lambda\in I_{x\otimes y}$. Then $\lambda \leq x\otimes y$. Therefore, for every $\alpha\in S_x, \beta\in S_y$, $\lambda\leq x\otimes y\leq \alpha\otimes \beta$. This shows $\ess(x\otimes y)=\sup I_{x\otimes y}\leq \inf (S_x\otimes S_y)=\ess x\otimes\ess y$.  Let us now take $\lambda\in S_{x\otimes y}$. Then, $x\otimes y\leq\lambda$. Then, for every $\alpha\in I_x,\beta\in I_y$, we have $\alpha\otimes\beta\leq x\otimes y\leq \lambda$. This shows $\ess x\otimes\ess y=\sup(I_x\otimes I_y)\leq  \inf S_{x\otimes y}=\ess(x\otimes y)$.
	
	The function ${\rm ess}$ is surjective, because it must preserve every contant symbol $\lambda\in [0,1]$.
	
	By Lemma~\ref{l:then essinf=esssup}, if $A$ is subdirectly irreducible, then, for all $x\in A$, $\essinf x=\esssup x$.
\endproof 

We call an algebra $A$ \emph{simple} if it is not trivial and any proper quotient of $A$ is trivial. From Theorem~\ref{t:morphism to [0,1]}, we obtain that $[0,1]$ is the unique simple MC-algebra, as stated in Corollary~\ref{c:simple} below. This is similar to H{\"o}lder's Theorem for lattice-ordered groups.
\corollary\label{c:simple}
	Let $A\in \MC$. $A$ is simple if, and only if, it is isomorphic to $[0,1]$.
\endcorollary
\proof 
	$[0,1]$ is simple. Indeed, if $B\in\MC$, and $\varphi\colon [0,1]\rightarrow B$ is a surjective not-injective MC-morphism, then there exist $\alpha,\beta\in [0,1]$ such that $\alpha\neq \beta$ and $\varphi(\alpha)=\varphi(\beta)$, i.e., $\alpha_B= \beta_B$. By Lemma~\ref{l:trivial if and only if 0=1}, $B$ is trivial.
	Hence, $[0,1]$ is simple, as well as any of its isomorphic copies. Suppose that $A$ is simple. By the subdirect representation theorem, it is isomorphic to a subdirect product $\prod_{i\in I}B_i$ of subdirectly irreducible algebras. Since $A$ is simple, it is not trivial. Hence, $I\neq \emptyset$, and thus there exists a surjective morphism $\varphi\colon A\rightarrow B$, with $B$ subdirectly irreducible. By Theorem~\ref{t:morphism to [0,1]}, we have a surjective MC-morphism ${\rm ess}\colon B\rightarrow [0,1]$. Hence, we have a surjective MC-morphism $\psi\colon A\rightarrow [0,1]$, which must be injective since $A$ is simple. Hence $\psi$ is an isomorphism.
\endproof 
Corollary~\ref{c:simple} implies that the set of morphisms from an MC-algebra $A$ to $[0,1]$ is in bijection with the set of maximal congruences on $A$. This explains why we gave the name $\Max(A)$ for the set $\hom_{\MC}(A,[0,1])$. 
\subsection{The map ${\rm ess}$ preserves distance}

The main goal of this subsection is to prove that, for $A$ a nontrivial subdirectly irreducible algebra, the map ${\rm ess}\colon A\rightarrow [0,1]$ preserves $\dist$, i.e., $\dist_A(a,b)=\dist_{[0,1]}(\ess a, \ess b)$. We will actually prove, in Lemma~\ref{l:du is preserved}, a slightly stronger statement, i.e., that $\du$ is preserved.

\lemma\label{l:equiv plus minus}
	Let $A\in\MC$ and $x,y\in A$. Then $y\leq x\oplus\lambda$ if, and only if, $y\ominus\lambda\leq x$.
\endlemma
\proof 
	If $y\leq x\oplus\lambda$, then $y\ominus\lambda\leq(x\oplus \lambda)\ominus \lambda\leq x$. If $y\ominus\lambda\leq x$, then 
	$y\leq (y\ominus \lambda)\oplus\lambda\leq x\oplus \lambda$.
\endproof 

\remark
	Let $A\in\MC$ and $x,y\in A$. Then, by Remark~\ref{l:equiv plus minus}, we have
	$$\uparrow_x^y\coloneqq \{\lambda\in[0,1]\mid y\leq x\oplus \lambda\}=\{\lambda\in[0,1]\mid y\ominus \lambda\leq x\}.$$ 
	Hence
	$$\du(x,y)\coloneqq \inf\uparrow_x^y= \inf\{\lambda\in[0,1]\mid y\leq x\oplus \lambda\}=\inf\{\lambda\in[0,1]\mid y\ominus \lambda\leq x\}.$$    
\endremark

\remark
	Let $A\in\MC$. Then the set $\uparrow_x^y$ calculated in $A$ equals $\uparrow_y^x$ calculated in $A^\partial$. Thus $\du_A(x,y)=\du_{A^\partial}(y,x)$.
\endremark

\remark
	\begin{enumerate}
		\item Let $I$ be a set, and, for each $i\in I$, let $A_i\in \MC$. Let $a,b\in \prod_{i\in I}A_i$. Then, $\du(a,b)=\sup_{i\in I}\du(a_i,b_i)$.
		\item Let $A\in\MC$, let $B$ be an MC-subalgebra of $A$, and let $x,y\in B$. Then $\du(x,y)$ is the same calculated in $A$ and $B$.
	\end{enumerate}	
\endremark 

\lemma
	Let $A\in\MC$ and $x,y,z\in A$. The following properties hold.
	\begin{enumerate}
		\item\label{i:in-interval} $\dist^{\uparrow}(x,y)\in [0,1]$.
		\item\label{i:uparrow-zero} $x\leq y\Rightarrow\du(y,x)=0$.
		\item\label{i:triang-ineq} $\dist^{\uparrow}(x,z)\leq \dist^{\uparrow}(x,y)\oplus\dist^{\uparrow}(y,z)\leq \dist^{\uparrow}(x,y)+\dist^{\uparrow}(y,z) $.
	\end{enumerate}
\endlemma
\proof 
	\eqref{i:in-interval} and \eqref{i:uparrow-zero} are clear, by definition. To prove \eqref{i:triang-ineq}, let $\lambda_0\in\uparrow_x^y$ and $\lambda_1\in \uparrow_y^z$. Then $y\leq x\oplus \lambda_0$, $z\leq y\oplus \lambda_1$. Then $z\leq y\oplus \lambda_1\leq  (x\oplus\lambda_0)\oplus\lambda_1=x\oplus(\lambda_0\oplus \lambda_1)$. Therefore, $\lambda_0\oplus \lambda_1\in \uparrow_x^z$. Therefore, $\du(x,z)=\inf(\uparrow_x^z)\leq\inf (\uparrow_x^y\oplus \uparrow_y^z)$. By continuity of $\oplus$, $\inf (\uparrow_x^y\oplus \uparrow_y^z)=\inf (\uparrow_x^y)\oplus \inf(\uparrow_y^z)=\du(x,y)\oplus \du(y,z)$.
\endproof

\lemma
	Let $A\in \MC$. Then
	\begin{enumerate}
		\item For every $x\in A$, $\dist^{\uparrow}(x,\essinf x)=\dist^{\uparrow}(\esssup x,x)=0$.
		
		\item If $A$ is nontrivial, for every $\alpha,\beta\in [0,1]$, $\dist^{\uparrow}(\alpha_A,\beta_A)=(\beta-\alpha)^+$. 
	\end{enumerate}
\endlemma
\proof 
	\begin{enumerate}		
		\item Let $\varepsilon>0$. Then $ x\leq \esssup x\oplus \varepsilon$. Therefore, $\varepsilon\in \uparrow_{\esssup x}^x$. Since it holds for every $\varepsilon$, then $\du(\esssup x,x)=\inf \uparrow_{\esssup x}^x=0$.
		Via the order-dual algebra, $\dist^{\uparrow}(x,\essinf x)=0$ is automatically proven.
		\item If $A$ is nontrivial, then $\uparrow_\alpha^\beta=\{\lambda\in[0,1]\mid \beta_A\leq \alpha_A\oplus \lambda\}=\{\lambda\in[0,1]\mid \lambda\geq (\beta-\alpha)^+ \}$.
	\end{enumerate}
\endproof

\lemma\label{l:du is preserved}
	Let $A\in \MC$, and let $x,y\in A$ be such that $\essinf x=\esssup x$ and $\essinf y=\esssup y$. Then, $\dist^{\uparrow}(x,y)=(\ess y -\ess x)^+$.
\endlemma
\proof 
	$A$ is nontrivial, because $\essinf x=\esssup x$. Thus, $(\ess y -\ess x)^+=\du(\ess x,\ess y)$. We are left to prove $\dist^{\uparrow}(x,y)=\du(\ess x,\ess y)$. We have
	$$\mld 
	\dist^{\uparrow}(x,y)&\leq \dist^{\uparrow}(x,\ess x)\oplus\dist^{\uparrow}(\ess x,\ess y)\oplus\dist^{\uparrow}(\ess y,y)=\\
	=0\oplus\dist^{\uparrow}(\ess x,\ess y)\oplus0=\\
	=\dist^{\uparrow}(\ess x,\ess y).
	$$
	Moreover, 
	$$\mld \dist^{\uparrow}(\ess x,\ess y)&\leq \du(\ess x,x)\oplus \du(x,y)\oplus \du(y,\ess y)=\\
	=0\oplus \du(x,y)\oplus 0=\\
	=\du(x,y).
	$$
\endproof 

\subsection{Every archimedean MC-algebra is an algebra of functions}
In this subsection, we prove that every archimedean MC-algebra has enough morphisms towards $[0,1]$ to separate its elements; this completes the proof of Theorem~\ref{t:arch iff}.

\lemma\label{l:morphism on d}
	Let $A,B\in \MC$ and let $\phi\colon A\rightarrow B$ be an MC-morphism. Then, for every $x,y\in A$,
	$$\du(\phi(x),\phi(y))\leq \dist^{\uparrow}(x,y).$$
\endlemma
\proof 
	Let $\lambda\in \uparrow_x^y$: $y\leq x\oplus \lambda$. Then $\phi(y)\leq \phi(x)\oplus \lambda$; thus $\lambda\in \uparrow_{\phi(x)}^{\phi(y)}$. Therefore, $\dist\uparrow(\phi(x),\phi(y))=\inf\uparrow_{\phi(x)}^{\phi(y)}\leq \inf\uparrow_x^y= \dist^{\uparrow}(x,y)$.
\endproof 

\theorem \label{t:d is sup morph}
	Let $A\in \MC$, and let $x,y\in A$. Then,
	$$\dist^{\uparrow}(x,y)=\sup_{\varphi\colon A\rightarrow [0,1] \ \text{MC-morphism}}( \varphi(y)-\varphi(x))^+.$$
\endtheorem
\proof 
	By the subdirect representation theorem, $A$ is an MC-subalgebra of a product of subdirectly irreducible MC-algebras. Say $\iota\colon A\rightarrow \prod_{i\in I}A_i$. For each $i\in I$, consider the projection $\pi_i\colon \prod_{i\in I}A_i\rightarrow A_i$ and the morphism $\ess_i\colon A_i\rightarrow [0,1]$ as in Theorem~\ref{t:morphism to [0,1]}.
	Then 
	$$\dist^{\uparrow}(x,y)=\sup_{i\in I}\dist^{\uparrow}(\pi_i\iota(x),\pi_i\iota(y))=\sup_{i\in I}(\ess_i\pi_i\iota(y)-\ess_i\pi_i\iota(x))^+.$$
	Since $\ess_i\circ \pi_i\circ \iota\colon A\rightarrow [0,1]$ is an MC-morphism, we obtain 
	$$\dist^{\uparrow}(x,y)\leq\sup_{\varphi\colon A\rightarrow [0,1] \ \text{MC-morphism}}(\varphi(y)-\varphi(x))^+.$$
	From Lemma~\ref{l:morphism on d}, we obtain the converse inequality.
\endproof

\lemma
	Let $A\in\MC$ and $x,y,z\in A$. The following properties hold.
	\begin{enumerate}
		\item\label{i:in-interval-dist} $\dist(x,y)\in [0,1]$.
		\item\label{i:antisym} $\dist(x,x)=0$.
		\item\label{i:triang-ineq-dist} $\dist(x,z)\leq \dist(x,y)\oplus\dist(y,z)\leq\dist(x,y)+\dist(y,z).$
	\end{enumerate}
\endlemma
\proof 
	\eqref{i:in-interval-dist} and \eqref{i:antisym} are clear. Let us prove \eqref{i:triang-ineq-dist}. 
	$$
	\mld \dist(x,z)&=\max\{\du(x,z),\du(z,x) \}\leq\\
	\leq\max\{\du(x,y)\oplus\du(y,z),\du(z,y)\oplus\du(y,x)\}\leq\\
	\leq\max\{\du(x,y), \du(y,x) \}\oplus\max\{\du(y,z), \du(z,y) \}=\\
	=\dist(x,y)\oplus\dist(y,z).
	$$
\endproof 
\remark
	Let $A\in\MC$. Then $A$ is archimedean if, and only if, $(A,\dist)$ is a metric space.
\endremark 
\theorem \label{t:dist as sup}
	Let $A\in \MC$, and let $x,y\in A$. Then
	$$\dist(x,y)=\sup_{\varphi\colon A\rightarrow [0,1] \ \text{MC-morphism}} \rvert\varphi(x)-\varphi(y)\rvert.$$
\endtheorem
\proof 
	$$
	\mld \dist(x,y)&=\max\{\du(x,y), \du(y,x) \}\stackrel{\text{Thm.~\ref{t:d is sup morph}}}{=}\\
	=\max\left\{ \sup_{\varphi\colon A\rightarrow [0,1] \ \text{MC-morph.}}( \varphi(y)-\varphi(x))^+,\sup_{\varphi\colon A\rightarrow [0,1] \ \text{MC-morph.}}( \varphi(x)-\varphi(y))^+ \right\}=\\
	=\sup_{\varphi\colon A\rightarrow [0,1] \ \text{MC-morph.}}\max\left\{( \varphi(y)-\varphi(x))^+,( \varphi(x)-\varphi(y))^+\right\}=\\
	=\sup_{\varphi\colon A\rightarrow [0,1] \ \text{MC-morph.}} \rvert\varphi(x)-\varphi(y)\rvert.
	$$
\endproof

We are ready to prove the implication [\eqref{i:arch}$\Rightarrow$\eqref{i:enough-morph}] in Theorem~\ref{t:arch iff}.
\theorem \label{t:arch iff separates points}
	Let $A\in \MC$ be archimedean, and let $x,y\in A$ with $x\neq y$. Then, there exists an MC-morphism $\varphi\colon A\rightarrow [0,1]$ such that $\varphi(x)\neq \varphi(y)$.
\endtheorem
\proof 
	By Theorem~\ref{t:dist as sup}, for every $x,y\in A$, we have
	$$\dist(x,y)=\sup_{\varphi\colon A\rightarrow [0,1] \ \text{MC-morphism}}\lvert \varphi(x)-\varphi(y)\rvert.$$
	Since $A$ is archimedean, $\dist(x,y)\neq 0$. Therefore, $\sup_{\varphi\colon A\rightarrow [0,1] \ \text{MC-morphism}}\lvert \varphi(x)-\varphi(y)\rvert\neq 0$, and hence there exists $\varphi\colon A\rightarrow [0,1]$ MC-morphism such that $\lvert\varphi(x)-\varphi(y)\rvert\neq 0$, i.e., $\varphi(x)\neq \varphi(y)$.
\endproof 

\corollary\label{c:ev preserves distances}
	For every $A\in \MC$, and every $x,y\in A$, we have $\dist_A(x,y)=\dist_{\C\Max(A)}(\ev_x,\ev_y)$.
\endcorollary

The results obtained so far enable us to prove one of the main results of this section---namely, Theorem~\ref{t:arch iff}.

\proof [of Theorem~\ref{t:arch iff}]
	By Remark~\ref{r:proof arch} and Theorem~\ref{t:arch iff separates points}.
\endproof

We add one additional characterization of archimedean MC-algebras.

\theorem 
	Let $A\in\MC$. Then the following conditions are equivalent.
	\begin{enumerate}
		\item\label{i:archim} $A$ is archimedean.
		\item\label{i:dist-leq} For every $x,y\in A$, $\du(x,y)=0$ implies $y\leq x$.
	\end{enumerate}
\endtheorem
\proof 
	\noindent {[\eqref{i:dist-leq}$\Rightarrow$\eqref{i:archim}]} Let $x,y\in A$, and suppose $\dist(x,y)=0$. Then $\du(x,y)=0$ and $\du(y,x)=0$. Hence $y\leq x$ and $x\leq y$. Therefore $x=y$.
	
	\noindent [\eqref{i:archim}$\Rightarrow$\eqref{i:dist-leq}] Let us suppose $A$ is archimedean, and let $x,y\in A$ be such that $\dist^{\uparrow}(x,y)=0$. Then, for every MC-morphism $\varphi\colon A\rightarrow [0,1]$, $\varphi(y)\leq \varphi(x)$; hence $\varphi(y)\lor \varphi(x)=\varphi(x)$, which implies $\varphi(y\lor x)=\varphi(x)$. Hence, for all MC-morphisms $\varphi\colon A\rightarrow [0,1]$, $\lvert \varphi(y\lor x)-\varphi(x)\rvert=0$. By Theorem~\ref{t:dist as sup}, $\dist(y\lor x,x)=0$. Since $A$ is archimedean, $y\lor x=x$, that is, $y\leq x$.
\endproof

\section{$A$ is Cauchy complete if, and only if, the unit $\varepsilon_A$ is surjective}

The aim of this section is to prove Theorem~\ref{t:Cauchy complete iff} above, which states, for any $A\in \MC$, the equivalence of the following conditions.
\begin{enumerate}
	\item\label{i:complete} $A$ is Cauchy complete.
	\item\label{i:surj} The unit $\varepsilon_A\colon A\rightarrow \C\Max(A)$ is surjective.
\end{enumerate}

\remark\label{r:surj then Cauchy}
	The implication [\eqref{i:surj}$\Rightarrow$\eqref{i:complete}] of Theorem~\ref{t:Cauchy complete iff}---i.e., if the unit $\varepsilon_A\colon A\rightarrow \C\Max(A)$ is surjective, then $A$ is Cauchy complete---follows from the fact, observed in Corollary~\ref{c:ev preserves distances}, that $\varepsilon_A$ preserves $\dist$. In detail, let $(a_n)_{n\in \N}$ be a Cauchy sequence in $A$. Then $(\varepsilon_A(a_n))_{n\in \N}$ is a Cauchy sequence in $\C{\Max{(A)}}$. Since $\C\Max(A)$ is Cauchy complete (see [\eqref{i:in-the-image-alg}$\Rightarrow$\eqref{i:arch-cc}] in Theorem~\ref{t:characterization algebraic fixed objects}), there exists $f\in \C\Max(A)$ such that $(\varepsilon_A(a_n))_{n\in \N}$ converges to $f$. Since $\varepsilon_A$ is surjective, there exists $a\in A$ such that $\varepsilon_A(a)=f$. The sequence $(a_n)_{n\in \N}$ converges to $a$.
\endremark

We are left to prove that if $A$ is Cauchy complete, then $\varepsilon_A$ is surjective. To do so, we make use of an analogue of Stone-Weierstrass Theorem. Our source of inspiration is \citep[Section 7]{HofmannLong}.
\lemma\label{l:separated} 
	Let $X$ be a preordered topological space, let $x,y\in X$, let $L$ be an MC-subalgebra of $\C(X)$, and let $\phi\in L$ be such that $\phi(x)<\phi(y)$. Then, there there exists $\psi\in L$ and an open neighbourhood $U_y$ of $y$ such that $\psi(x)=0$ and, for all $z\in U_y$, $\psi(z)=1$.
\endlemma
\proof 
	There exists $c\in [0,1]$ such that $\varphi(x)<c<\varphi(y)$. Let $n\in \N$ be such that $n(c-\varphi(x))\geq 1$. Set $\psi\coloneqq\underbrace{(\varphi\ominus\varphi(x))\oplus \dots\oplus (\varphi\ominus\varphi(x))}_{n \text{ times}}$. Set $U_y$ as the pre-image of $(c,1]$ under $\varphi$. Note that $y\in U_y$. We have $\psi(x)=\underbrace{(\varphi(x)\ominus\varphi(x))\oplus \dots\oplus (\varphi(x)\ominus\varphi(x))}_{n \text{ times}}=0$, and, for every $z\in U_y$, $\psi(z)=\underbrace{(\varphi(z)\ominus\varphi(x))\oplus \dots\oplus (\varphi(z)\ominus\varphi(x))}_{n \text{ times}}=1$.
\endproof 
\theorem[Ordered version of Stone-Weierstrass Theorem]\label{t:StWe}
	Let $X$ be a preordered topological space, let $L$ be an MC-subalgebra of $\C(X)$, and suppose that, for every $x,y\in X$, if $x\not\geq y$ then there exists $\phi\in L$ such that $\phi(x)<\phi(y)$. If $X$ is compact, then, for every $\psi\in \C(X)$, there exists a sequence $(\phi_n)_{n\in \N}$ in $L$ converging to $\psi$ in the sup metric.
\endtheorem
\proof 
	Fix $\varepsilon\in (0,1]$; we shall find $\phi\in L$ such that $\sup_{x\in X}\lvert \psi(x)-\phi(x)\rvert\leq \varepsilon$. Fix $x\in X$. Set $U\coloneqq\{z\in X\mid \psi(z)<\psi(x)+ \varepsilon \}$. The set $U$ is open. Moreover, for every $y\in X$ such that $y\leq x$, we have $y\in U$ (by monotonicity of $\psi$); contrapositively, for every $y\in X\setminus U$ we have $x\not\geq y$. Hence, by Lemma \ref{l:separated}, for every $y\in X\setminus U$ there exists $\alpha_y\in L$ and an open neighbourhood $U_y$ of $y$ such that $\alpha_y(x)=0$ and, for all $z\in U_y$, $\alpha_y(z)=1$.
	
	By compactness of $X$, there exist finitely many elements $y_1,\dots,y_n\in X\setminus U$ such that $X=U\cup U_{y_1}\cup \dots\cup U_{y_n}$. Set $\lambda\coloneqq \psi(x)$, and set $\overline{\lambda}\colon X\rightarrow [0,1]$ to be the function constantly equal to $\lambda$. Let us define $\phi_x\coloneqq \alpha_{y_1}\oplus \dots\oplus \alpha_{y_n}\oplus \overline{\lambda}$.
	We claim that $\phi_x$ has the following properties.
	\begin{enumerate}[label=a\arabic*., ref=a\arabic*]
		\item\label{i:psi-in-x} $\phi_x(x)=\psi(x)$.
		\item\label{i:approx} For every $z\in X$, $\phi_x(z)>\psi(z)- \varepsilon$.
	\end{enumerate}
	Indeed, \eqref{i:psi-in-x} holds because, for $1\leq i\leq n$, we have $\alpha_{y_i}(x)=0$, and so $\phi_x(x)= \alpha_{y_1}(x)\oplus \dots\oplus \alpha_{y_n}(x)\oplus \lambda=0\oplus \dots\oplus 0\oplus \lambda=\lambda=\psi(x)$. We prove \eqref{i:approx} by cases. If $z\in U$, then $\psi(x)=\alpha_{y_1}(z)\oplus \dots\oplus \alpha_{y_n}(z)\oplus\lambda\geq \lambda=\psi(z)>\psi(z)-\varepsilon$. If $z\in X\setminus U$, there exists $i\in \{1,\dots,n \}$ such that $z\in U_{y_i}$. Thus, $\phi_x(z)=\alpha_{y_1}(z)\oplus \dots\oplus \alpha_{y_n}(z)\oplus \lambda=1\geq \psi(z)>1-\varepsilon$. This settles the claim that \eqref{i:psi-in-x} and \eqref{i:approx} hold.

	Now $x$ is not fixed anymore. For $x\in X$, set
	$$V_x\coloneqq\{z\in X\mid \phi_x(z)<\psi(z)+\varepsilon \}.$$
	The set $V_x$ is open because the functions $\phi_x$ and $\psi$ are continuous. Moreover $x\in V_x$ because of \eqref{i:psi-in-x}. Therefore the family $(V_x)_{x\in X}$ is an open cover of $X$. Again, by compactness of $X$, there exists a finite subcover $V_{x_1}, \dots, V_{x_m}$ of $X$. Define $\phi\coloneqq\phi_{x_1}\land\dots \land \phi_{x_m}$; note that $\phi\in L$. For all $z\in X$ we have the following.
	\begin{enumerate}[label=b\arabic*., ref=b\arabic*]
		\item There exists $i\in \{1,\dots, n \}$ such that $z\in V_{x_i}$. Hence,  $\phi(z)=\phi_{x_1}(z)\land\dots \land \phi_{x_m}(z)\leq \phi_{x_i}(z)<\psi(z)+\varepsilon.$
		\item By \eqref{i:approx}, $\phi(z)=\phi_{x_1}(z)\land\dots\land\phi_{x_m}(z)>(\psi(z)-\varepsilon)\land\dots\land(\psi(z)-\varepsilon)=\psi(z)-\varepsilon$.
	\end{enumerate}
	Hence, for all $z\in X$, $\psi(z)-\varepsilon<\phi(z)<\psi(z)+\varepsilon$, which implies $\sup_{x\in A}\lvert \psi(x)-\phi(z)\rvert\leq \varepsilon$.
	
\endproof 

\theorem \label{t:epsilon is surj}
	Let $A\in\MC$. If $A$ is Cauchy complete, then 
	\begin{align*}
	\varepsilon_A\colon A&\longrightarrow \C\Max(A)\\
	a&\longmapsto \ev_a\colon \Max(A)\rightarrow [0,1],\ x\mapsto x(a)
	\end{align*}
	is surjective.
\endtheorem
\proof 
	Set $L\subseteq \C\Max(A)$ as the image of $A$ under $\varepsilon_A$. $L$ is an MC-subalgebra of $\C\Max(A)$. By the definition of the topology and the order on $\Max(A)$, the hypothesis in Theorem~\ref{t:StWe} are fulfilled, with $X\coloneqq\Max(A)$. Hence Theorem~\ref{t:StWe} applies: for every $\psi\in \C\Max(A)$, there exists a sequence $(\widetilde{a}_n)_{n\in \N}$  in the image of $\varepsilon_A$ converging to $\psi$  with respect to the sup metric. Let $(a_n)_{n\in \N}$ be a sequence in $A$ such that, for every $n\in \N$, $\ev_{a_n}=\widetilde{a}_n$. Therefore, $(\widetilde{a}_n)$ is a Cauchy sequence with respect to the sup metric. Therefore, for every $\varepsilon>0 $, there exists $k\in\N$ such that, for every $n,m\geq k $, we have $\dist(\ev_{a_{n}},\ev_{a_{m}})<\varepsilon$. Since, by Corollary~\ref{c:ev preserves distances}, $\varepsilon_A$ preserves $\dist$, $\dist(a_n,a_m)=\dist(\ev_{a_{n}},\ev_{a_{m}})<\varepsilon$, and therefore $(a_n)_{n\in \N}$ is a Cauchy sequence. Since $A$ is Cauchy complete, there exists $a\in A$ such that $a_n$ converges to $a$. Therefore, for every $\varepsilon>0$, there exists $n\in \N$ such that, for all $m\geq n$, $\dist(a_n,a)<\varepsilon$. Since, by Corollary~\ref{c:ev preserves distances}, $\varepsilon_A$ preserves $\dist$,  for every $\varepsilon>0$, there exists $n\in \N$ such that, for all $m\geq n$, $\dist(\ev_{a_n},\ev_a)=\dist(a_n,a)<\varepsilon$. Hence, $\ev_{a_n}$ converges both to $\ev_a$ and $\psi$. Therefore, $\psi=\ev_a$.
\endproof 

We can now prove Theorem~\ref{t:Cauchy complete iff}.
\proof [of Theorem~\ref{t:Cauchy complete iff}]
	By Remark~\ref{r:surj then Cauchy} and Theorem~\ref{t:epsilon is surj}.
\endproof 
Finally, we can prove Theorem~\ref{t:characterization algebraic fixed objects}.
\proof [of Theorem~\ref{t:characterization algebraic fixed objects}]
	By Remark~\ref{r:easy} and Theorems~\ref{t:arch iff} and~\ref{t:Cauchy complete iff}.
\endproof 
Combining Theorems~\ref{t:characterization fixed geometric} and~\ref{t:characterization algebraic fixed objects}, we have that the adjoint contravariant functors $\C\colon \PreT\rightarrow \MC$ and $\Max\colon \MC\rightarrow \PreT$ restrict to a dual equivalence between $\PosComp$ and the full subcategory of archimedean Cauchy complete MC-algebras. Hence, we have the following.
\theorem \label{t:equivalence of cat}
	The dual of $\PosComp$ is equivalent to the full subcategory of $\MC$ given by the archimedean Cauchy complete MC-algebras.
\endtheorem

\section{The variety $\MCinfty$}

Up to now, we have proved that the category $\PosComp$ of compact ordered spaces is dually equivalent to the full subcategory of archimedean Cauchy complete MC-algebras. Our final goal is to show that the full subcategory  of archimedean Cauchy complete MC-algebras is isomorphic to a variety. Our strategy to achieve this purpose is analogous to the strategy that, in \citep{HofmannShort}, Section 3, after Theorem 3.8, was pursued to show that a given category was a quasi-variety. The crucial difference is that we make the axioms equational; in achieving the equational axiomatization, having the operation $\land$ amongst the primitive operations has facilitated us. 

\subsection{Adding the infinitary ``Cauchy'' operation $\delta$}
In order to ensure Cauchy completeness, we would like to add an operation $\delta$ of countably infinite arity to the class of operations of $\MC$ that computes the limit of ``enough'' Cauchy sequences, meaning that convergence of such sequences in an MC-algebra is enough to imply Cauchy completeness (and, at the same time, it is possible to interpret $\delta$ in $[0,1]$ so that it becomes a monotone continuous function from $[0,1]^\N$ to $[0,1]$ that calculates the limit of such sequences, see \citep[p. 283]{HofmannShort}).
\definition
	Let $A\in\MC$. A sequence $(a_n)_{n\in \N}$ in $A$ is called \emph{HNN-Cauchy} if, for every $n\in \N$, 
	$$a_n\leq a_{n+1}\leq a_n\oplus\frac{1}{2^n}.$$
\enddefinition
\noindent This definition is inspired by Lemma 3.9 in \citep{HofmannShort}; in fact, ``HNN'' stands for ``Hofmann, Neves, Nora'', the authors of the paper.

\lemma\label{l:supercauchy is cauchy}
	Let $A\in \MC$, and let $(a_n)_{n\in \N}$ be an HNN-Cauchy sequence in $A$. Then, for every $n,m\in \N$, with $n\leq m$, we have
	$$a_n\leq a_m\leq a_n\oplus\frac{1}{2^{n-1}},$$
	and therefore $(a_n)_{n\in \N}$ is a Cauchy sequence.
\endlemma
\proof 
	The inequality $a_n\leq a_m$ is obtained by induction on $m$. Moreover, 
	$$a_m\leq  (\dots(a_{n}\oplus \frac{1}{2^n})\oplus\dots)\oplus \frac{1}{2^{m-1}}\leq a_n\oplus \sum_{i=n}^\infty\frac{1}{2^i}=a_n\oplus \frac{1}{2^{n-1}}.$$
	Hence, 
	$$\dist(a_n,a_m)=\max\{\du(a_n,a_m), \du(a_m,a_n)\}=\max\{\du(a_n,a_m),0 \}=\du(a_n,a_m)\leq\frac{1}{2^{n-1}}.$$
\endproof 

\lemma
	For $A\in \MC$, the following conditions are equivalent.
	\begin{enumerate}
		\item\label{i:cauchy-complete} $A$ is Cauchy complete.
		\item\label{i:converges} Every HNN-Cauchy sequence in $A$ converges.
	\end{enumerate}
\endlemma
\proof 
	{	[\eqref{i:cauchy-complete}$\Rightarrow$\eqref{i:converges}]} By Lemma~\ref{l:supercauchy is cauchy}.
	
	\noindent [\eqref{i:converges}$\Rightarrow$\eqref{i:cauchy-complete}] Let $(a_n)_{n\in \N}$ be a Cauchy sequence. For every $i\in \N$, let $k_i\in \N$ be such that, for every $n,m\geq k_i$, $\dist( a_n,a_m)< \frac{1}{2^{i+1}}$. For each $i\in \N$, set $b_i\coloneqq a_{k_i}$. Then, for every $i\in \N$, and every $n,m\geq i$, $\dist( b_n,b_m)< \frac{1}{2^{i+1}}$. In particular, for each $i\in \N$, $\dist(b_{i}, b_{i+1})< \frac{1}{2^{i+1}}$. Therefore,
	$$b_{i}\ominus \frac{1}{2^{i+1}}\leq b_{i+1};$$
	$$b_{i+1}\ominus \frac{1}{2^{i+1}}\leq b_{i}.$$
	For all $i\in \N$, we set $c_i\coloneqq b_i\ominus \frac{1}{2^i}$.
	
	Then, for every $i\in \N$,  
	$$c_{i}=b_{i}\ominus \frac{1}{2^{i}}=\left(b_{i}\ominus \frac{1}{2^{i+1}}\right)\ominus \frac{1}{2^{i+1}}\leq b_{i+1}\ominus \frac{1}{2^{i+1}}=c_{i+1};$$
	$$c_{i+1}=b_{i+1}\ominus \frac{1}{2^{i+1}}\leq b_{i}\leq \left(b_{i}\ominus \frac{1}{2^{i}}\right)\oplus \frac{1}{2^{i}}=c_{i}\oplus \frac{1}{2^{i}}.$$
	So, $c_{i}\leq c_{i+1}\leq c_{i}\oplus \frac{1}{2^{i}}$. Hence, $(c_n)_{n\in \N}$ is an HNN-Cauchy sequence, and thus there exists $c\in A$ such that $(c_n)_{n\in \N}$ converges to $c$. 
	
	We have 
	$$\dist(b_i,c)\leq \dist(b_i,c_i)+\dist(c_i,c)=\dist\left(b_i,b_i\ominus \frac{1}{2^i}\right)+\dist(c_i,c)\leq \frac{1}{2^i}+\dist(c_i,c)\stackrel{i\rightarrow \infty}{\rightarrow }0.$$
	Therefore, the sequence $(b_i)_{n\in \N}$ converges to $c$. The sequence $(a_n)_{n\in \N}$ is a Cauchy sequence that admits a convergent subsequence $(b_i)_{i\in \N}$; by a standard argument, it follows that $(a_n)_{n\in \N}$ converges.
\endproof 

\notation
	Inductively on $n\in \N$, we define the term $\rho_n$ of arity $n+1$ in the language of $\MC$ as follows.
	$$\rho_0(x_0)\coloneqq x_0;$$
	$$\text{for $n\in \N$}\ \ \rho_{n+1}(x_0,\dots,x_{n+1})\coloneqq (x_0\lor \dots \lor x_{n+1})\land \left(\rho_{n}(x_0,\dots,x_{n})\oplus \frac{1}{2^{n}}\right).$$
\endnotation

\lemma\label{l:properties of rho}
	Let $A\in \MC$. For every $n\in \N$, and every $x_0,\dots,x_{n+1}\in A$, the following properties hold.
	$$\rho_{n}(x_0,\dots,x_{n})\leq \rho_{n+1}(x_0,\dots,x_{n+1})\leq \rho_{n}(x_0,\dots,x_{n})\oplus\frac{1}{2^{n}} .$$	
\endlemma
\proof 
	By definition of $\rho_n$, we have $\rho_{n}(x_0,\dots,x_{n})\leq x_0\lor \dots\lor x_{n}$ and $\rho_{n+1}(x_0,\dots,x_{n+1})\leq \rho_{n}(x_0,\dots,x_{n})\oplus\frac{1}{2^{n}}$. As a consequence, $\rho_{n}(x_0,\dots,x_{n})\leq x_0\lor \dots\lor x_{n}\leq x_0\lor \dots\lor x_{n+1}$ and $\rho_{n}(x_0,\dots,x_{n})\leq \rho_{n}(x_0,\dots,x_{n})\oplus \frac{1}{2^{n}}$. Thus, $\rho_{n}(x_0,\dots,x_{n})\leq (x_0\lor \dots\lor x_{n+1})\land \left(\rho_{n}(x_0,\dots,x_{n})\oplus \frac{1}{2^{n}}\right)=\rho_{n+1}(x_0,\dots,x_{n+1})$.
\endproof 

\lemma\label{l:in V dist rho_n rho_m}
	Let $A\in \MC$, and let $(x_n)_{n\in \N}$ be a sequence in $A$. The following properties hold.
	\begin{enumerate}
		\item\label{i:HNN} The sequence $(\rho_n(x_0,\dots,x_n))_{n\in \N}$ is an HNN-Cauchy sequence.
		\item \label{i:projection} If $(x_n)_{n\in \N}$ is an HNN-Cauchy sequence, then, for all $n\in \N$, 
		$$\rho_n(x_0,\dots,x_n)=x_n.$$
	\end{enumerate}
\endlemma
\proof 
	\eqref{i:HNN} follows from Lemma~\ref{l:properties of rho}. \eqref{i:projection} is proved inductively. The case $n=0$ is trivial. Inductive step: let $n\in \N$. Then $\rho_{n+1}(x_0,\dots,x_{n+1})\coloneqq (x_0\lor \dots \lor x_{n+1})\land \left(\rho_{n}(x_0,\dots,x_{n})\oplus \frac{1}{2^{n}}\right)\stackrel{\text{ind. hyp.}}{=} (x_0\lor \dots \lor x_{n+1})\land \left(x_n\oplus \frac{1}{2^{n}}\right)=x_{n+1}\land \left(x_n\oplus \frac{1}{2^n}\right)=x_{n+1}$.
\endproof 

Let $A\in \MC$, and let $(a_n)_{n\in \N}$ be a sequence in $A$. If $A$ is Cauchy complete, the sequence $(\rho_n(a_0,\dots,a_n))_{n\in \N}$ admits a limit in $A$. If, additionally, $A$ is archimedean, $(A,\dist)$ is a metric space and thus the limit is unique.
\notation
	Let $A\in \MC$ be archimedean and Cauchy complete. For every sequence $(a_n)_{n\in \N}$ in $A$, we set
	$$\delta(a_0,a_1,a_2,\dots)\coloneqq\lim_{n\rightarrow \infty} \rho_n(a_0,\dots,a_{n}).$$
\endnotation

The definition of $\delta$ takes inspiration from \citep[Section 3]{HofmannShort}; in fact, the function $\delta\colon [0,1]^\N\rightarrow [0,1]$ in \citep[Section 3]{HofmannShort} coincides with the interpretation in $[0,1]$ of what we call $\delta$ here.

\remark
	$\delta$ calculates the limit of HNN-Cauchy sequences.
\endremark  

\proposition\label{p:properties of W}
	Let $A\in \MC$ be archimedean and Cauchy complete. The following properties hold.
	\begin{enumerate}
		\item\label{i:constant} $\delta(x,x,x,\dots)=x$.
		\item\label{i:delta-lor} $\delta(x_0,x_1,x_2,\dots)\leq\delta(x_0\lor y_0,x_1\lor y_1, x_2\lor y_2,\dots)$.
		\item\label{i:delta-ominus} $\delta\left(x\ominus\frac{1}{2^0},x\ominus\frac{1}{2^1},x\ominus\frac{1}{2^2},\dots\right)=x$.
		\item\label{i:delta-sandwich} For all $n\in\N$ $$\rho_n(x_0,\dots,x_{n})\leq\delta(x_0,x_1,x_2,\dots)\leq\rho_n(x_0,\dots,x_{n})\oplus \frac{1}{2^{n-1}}.$$
	\end{enumerate}
\endproposition
\proof 
	Since $A\in \MC$ is archimedean, $A$ is (isomorphic to) a subalgebra of $[0,1]^X$, for some set $X$. The function $\dist\colon A\times A\rightarrow [0,1]$ coincides with the sup metric. We recall that, for every sequence $(f_n)_{n\in \N}$ in $A$, $\delta((f_n)_{n\in \N})=\lim_{n\rightarrow \infty} \rho_n(f_0,\dots,f_{n})$. The convergence is uniform, and therefore pointwise. Hence, it is enough to prove \eqref{i:constant}, \eqref{i:delta-lor}, \eqref{i:delta-ominus} and \eqref{i:delta-sandwich} for $A=[0,1]$.
	\begin{enumerate}
		\item The sequence $(x,x,x,\dots)$ is HNN-Cauchy; thus $\delta(x,x,x,\dots)=\lim_{n\rightarrow \infty} x=x$.
		\item For each $n\in \N$, set $z_n\coloneqq x_n\lor y_n$. By induction, we show $\rho_n(x_0,\dots,x_{n})\leq \rho_n(z_0,\dots,z_{n})$. Indeed, for $n=0$, we have $\rho_0(x_0)=x_0\leq z_0=\rho_0(z_0)$. Inductive step: let $n\in \N$; then, $\rho_{n+1}(x_0,\dots,x_{n+1})= (x_0\lor \dots \lor x_{n+1})\land \left(\rho_{n}(x_0,\dots,x_{n})\oplus \frac{1}{2^{n}}\right)\stackrel{\text{ind. hyp.}}{\leq} (z_0\lor \dots \lor z_{n+1})\land \left(\rho_{n}(z_0,\dots,z_{n})\oplus \frac{1}{2^{n}}\right)=\rho_{n+1}(z_0,\dots,z_{n+1})$. Hence, we have proved inductively $\rho_n(x_0,\dots,x_{n})\leq \rho_n(z_0,\dots,z_{n})$. Since, in $[0,1]$, $\lim$ is monotone, we have
		$$\mld
		\delta(x_0,x_1,x_2,\dots)&=\lim_{n\rightarrow \infty}\rho_n(x_0,\dots,x_n)\leq\\
		\leq  \lim_{n\rightarrow \infty}\rho_n(z_0,\dots,z_n)=\\
		=\delta(x_0\lor y_0,x_1\lor y_1, x_2\lor y_2,\dots).
		$$
		\item 
		Let us prove that $\left(x\ominus\frac{1}{2^0},x\ominus\frac{1}{2^1},x\ominus\frac{1}{2^2},\dots\right)$ is an HNN-Cauchy sequence. Indeed, $x\ominus\frac{1}{2^{n}}\leq x\ominus\frac{1}{2^{n+1}}\leq x\leq (x\ominus\frac{1}{2^{n}})\oplus \frac{1}{2^n}$.

		\item By Lemma~\ref{l:in V dist rho_n rho_m}, the sequence $(\rho_n(x_0,\dots,x_n))_{n\in \N}$ is an HNN-Cauchy sequence. By Lemma~\ref{l:supercauchy is cauchy}, for every $n,m\in \N$, with $n\leq m$, we have
		$$\rho_n(x_0,\dots,x_n)\leq \rho_m(x_0,\dots,x_m)\leq \rho_n(x_0,\dots,x_n)\oplus\frac{1}{2^{n-1}}.$$
		Fix $n$ and let $m$ tend to $\infty$.
	\end{enumerate}
\endproof

\bigskip

\subsection{The variety $\MCinfty$}

Recall the inductive definition of the term $\rho_n$ of arity $n+1$ in the language of $\MC$:
$$\rho_0(x_0)\coloneqq x_0;$$
$$\text{for $n\in \N$}\ \ \rho_{n+1}(x_0,\dots,x_{n+1})\coloneqq (x_0\lor \dots \lor x_{n+1})\land \left(\rho_{n}(x_0,\dots,x_{n})\oplus \frac{1}{2^{n}}\right).$$
\definition\label{d:Vd}
	We define the variety $\MCinfty$ as the variety obtained from the variety $\MC$ by adding an operation $\delta$ of countably infinite arity, together with the following additional axioms.
	\begin{enumerate}
		\item\label{i:constant-ax} $\delta(x,x,x,\dots)=x$.
		\item \label{i:delta-lor-ax} $\delta(x_0,x_1,x_2,\dots)\leq\delta(x_0\lor y_0,x_1\lor y_1, x_2\lor y_2,\dots)$.
		\item\label{i:delta-ominus-ax} $\delta\left(x\ominus\frac{1}{2^0},x\ominus\frac{1}{2^1},x\ominus\frac{1}{2^2},\dots\right)=x$.
		\item\label{i:delta-sandwich-ax} For all $n\in\N$ $$\rho_n(x_0,\dots,x_{n})\leq\delta(x_0,x_1,x_2,\dots)\leq\rho_n(x_0,\dots,x_{n})\oplus \frac{1}{2^{n-1}}.$$
	\end{enumerate}
\enddefinition

The idea behind this definition is the following. Axioms \eqref{i:constant-ax}, \eqref{i:delta-lor-ax}, \eqref{i:delta-ominus-ax} imply being archimedean (see Proposition~\ref{p:W is archimedean} below); Axiom \eqref{i:delta-sandwich-ax}  forces $\delta(x_0,x_1,x_2,\dots)$ to be the limit of $(\rho_n(x_0,\dots,x_n))_{n\in \N}$, and therefore it implies Cauchy completeness (see Proposition~\ref{p:W is complete} below).

\subsection{General properties of the forgetful functor $\MCinfty\rightarrow \MC$}
\proposition\label{p:W is archimedean}
	Let $A\in\MCinfty$. Then $A$ is archimedean.
\endproposition
\proof 
	Let $x,y\in A$ be such that $\dist(x,y)=0$. We shall show $x=y$. Since $\dist(x,y)=\max\{\du(x,y),\du(y,x) \}$, $\du(y,x)=0$. We recall $\du(y,x)=\inf \{\lambda\in[0,1]\mid x\ominus \lambda\leq y\}=0$. Hence, for all $\lambda\in(0,1]$, we have $x\ominus \lambda\leq y$.	Note that \eqref{i:delta-lor-ax} in Definition~\ref{d:Vd} implies that $\delta$ is monotone. We have $$x\stackrel{\eqref{i:delta-ominus-ax}}{=}\delta\left(x\ominus\frac{1}{2^0},x\ominus\frac{1}{2^1},x\ominus\frac{1}{2^2},\dots\right)\stackrel{\eqref{i:delta-lor-ax}}{\leq} \delta(y,y,y,\dots)\stackrel{\eqref{i:constant-ax}}{=}y.$$
	Analogously, one shows that $y\leq x$. Hence, $x=y$.
\endproof 

\theorem \label{t:forgetful is full}
	Let $A,B\in\MCinfty$, and let $\varphi\colon A\rightarrow B$ be an MC-morphism. Then, $\varphi$ preserves $\delta$.
\endtheorem
\proof 
	We should prove $\varphi(\delta_A(x_0,x_1,x_2,\dots))=\delta_B(\varphi(x_0),\varphi(x_1),\varphi(x_2),\dots)$. Since $B$ is archimedean, it is enough to prove $$\dist(\varphi(\delta_A(x_0,x_1,x_2,\dots)),\delta_B(\varphi(x_0),\varphi(x_1),\varphi(x_2),\dots))=0.$$
	For all $n\in\N$, we have 
	$$\rho_n(\varphi(x_0),\dots,\varphi(x_n))\leq \varphi(\delta_A(x_0,x_1,x_2,\dots))\leq \rho_n(\varphi(x_0),\dots,\varphi(x_n))\oplus \frac{1}{2^{n-1}},$$
	because $\varphi$ is an MC-morphism.	Moreover, since $B$ is an MC$_\infty$-algebra, we have
	$$\rho_n(\varphi(x_0),\dots,\varphi(x_n))\leq \delta_B(\varphi(x_0),\varphi(x_1),\varphi(x_2),\dots)\leq \rho_n(\varphi(x_0),\dots,\varphi(x_n))\oplus \frac{1}{2^{n-1}}.$$
	Hence, for all $n\in\N$, we have 
	$$\varphi(\delta_A(x_0,x_1,x_2,\dots))\leq \delta_B(\varphi(x_0),\varphi(x_1),\varphi(x_2),\dots)\oplus\frac{1}{2^{n-1}}$$
	and
	$$\delta_B(\varphi(x_0),\varphi(x_1),\varphi(x_2),\dots))\leq\varphi(\delta_A(x_0,x_1,x_2,\dots))\oplus\frac{1}{2^{n-1}}.$$
	Thus, 
	$$\dist(\varphi(\delta_A(x_0,x_1,x_2,\dots)),\delta_B(\varphi(x_0),\varphi(x_1),\varphi(x_2),\dots))=0.$$
\endproof

\proposition\label{p:W is complete}
	Let $A\in \MCinfty$. Then $A$ is Cauchy complete.
\endproposition
\proof 
	It is enough to prove that every HNN-Cauchy sequence in $A$ converges. Let $(x_n)_{n\in\N}$ be an HNN-Cauchy sequence in $A$. Then, $\rho_n(x_0,\dots,x_n)=x_n$. Hence, for all $n\in\N$, $x_n\leq \delta(x_0,x_1,x_2,\dots)\leq x_n\oplus \frac{1}{2^{n-1}}$, which implies $\dist(\delta(x_0,x_1,x_2,\dots), x_n)\leq \frac{1}{2^{n-1}}$, which implies that $\delta(x_0,x_1,x_2,\dots)$ is a limit for $(x_n)_{n\in\N}$.
\endproof

We denote with $\mathrm{U}_{\MCinfty,\MC}\colon \MCinfty\rightarrow \MC$ the forgetful functor.
\theorem \label{t:yesyesyes}
	\begin{enumerate}
		\item\label{i:full-faith} $\mathrm{U}_{\MCinfty,\MC}$ is full and faithful.
		\item\label{i:inj-on-obj} $\mathrm{U}_{\MCinfty,\MC}$ is injective on objects: $\mathrm{U}_{\MCinfty,\MC}(A)=\mathrm{U}_{\MCinfty,\MC}(B)$ implies $A=B$. This means that every MC-algebra admits at most one MC$_\infty$-structure that extends its MC-structure.
		\item\label{i:arch-and-Cauchy-complete} For $A\in\MC$, there exists $\widetilde{A}\in \MCinfty$ such that $\mathrm{U}_{\MCinfty,\MC}(\widetilde{A})=A$ if, and only if, $A$ is archimedean and Cauchy complete.
		\item\label{i:closed-under-iso} The image of $\mathrm{U}_{\MCinfty,\MC}$ on objects is closed under isomorphisms.
		\item\label{i:MCinfty} The MC-algebra $[0,1]$ admits a (unique) MC$_\infty$-structure.
		
	\end{enumerate}
\endtheorem
\proof 
	\begin{enumerate}
		\item The fact that $\mathrm{U}_{\MCinfty,\MC}$ is faithful is trivial. The fact that $\mathrm{U}_{\MCinfty,\MC}$ is full is Theorem~\ref{t:forgetful is full}.
		\item  Suppose $\mathrm{U}_{\MCinfty,\MC}(A)=\mathrm{U}_{\MCinfty,\MC}(B)$. Then, $A$ and $B$ share the same underlying set. Let $\mathrm{Id}\colon \mathrm{U}_{\MCinfty,\MC}(A)\rightarrow \mathrm{U}_{\MCinfty,\MC}(B)$ be the identity function (which is an MC-morphism since $\mathrm{U}_{\MCinfty,\MC}(A)=\mathrm{U}_{\MCinfty,\MC}(B)$). Since $\mathrm{U}_{\MCinfty,\MC}$ is full, there exists an MC$_\infty$-morphism $\varphi\colon A\rightarrow B$ such that $\mathrm{U}_{\MCinfty,\MC}(\varphi)=\mathrm{Id}$. Then, $\varphi$ is the identity function, and thus $A=B$. 
		\item Suppose there exists $\widetilde{A}\in \MCinfty$ such that $\mathrm{U}_{\MCinfty,\MC}(\widetilde{A})=A$. By Propositions~\ref{p:W is archimedean} and~\ref{p:W is complete}, $A$ is archimedean and Cauchy complete. For the converse implication, suppose that $A$ is archimedean and Cauchy complete. Then, by Proposition~\ref{p:properties of W}, $A$ admits an MC$_\infty$-structure.
		\item It is a consequence of \eqref{i:arch-and-Cauchy-complete}.
		\item $[0,1]$ is archimedean and Cauchy complete.
	\end{enumerate}
\endproof 
\corollary\label{c:isomorphism of categories}
	The variety $\MCinfty$ is isomorphic to the full subcategory of $\MC$ given by the archimedean Cauchy complete MC-algebras.
\endcorollary

We can now prove Theorem~\ref{t:MAIN-PosComp-is-variety}, which is our main result: the dual of $\PosComp$ is equivalent to a variety of algebras.
\proof [of Theorem~\ref{t:MAIN-PosComp-is-variety}]
	The dual of $\PosComp$ is equivalent, by Theorem~\ref{t:equivalence of cat}, to the full subcategory of $\MC$ given by the archimedean Cauchy complete MC-algebras, which is equivalent, by Corollary~\ref{c:isomorphism of categories}, to the variety $\MCinfty$. Hence, the dual of $\PosComp$ is equivalent to the variety $\MCinfty$.
\endproof

\section{The variety $\MCinfty$ and Linton's varietal theories}
\lemma\label{l:delta-cont-mon}
	The function $\delta\colon [0,1]^\N\rightarrow [0,1]$ is monotone and continuous with respect to the product order and product topology.
\endlemma
\proof 
	For every $n\in \N$, we set
	\begin{align*}
	\widetilde{\rho}_n\colon [0,1]^\N&\longrightarrow [0,1]\\
	(x_n)_{n\in \N}&\longmapsto \rho_n(x_0,\dots,x_{n}).
	\end{align*}
	Then, the sequence $(\widetilde{\rho}_n)_{n\in \N}$ tends to $\delta$ with respect to the the supremum norm, i.e., uniformly.	For every $i\in \N$, the projection $\pi_i\colon [0,1]^\N\rightarrow [0,1]$ onto the $i$-th coordinate is continuous, and for every $n\in \N$, $\rho_n\colon [0,1]^n\rightarrow [0,1]$ is continuous. Therefore, for every $n\in \N$, $\widetilde{\rho}_n\colon [0,1]^\N\rightarrow [0,1]$ is continuous. Since $(\widetilde{\rho}_n)_{n\in \N}$ converges to $\delta$ uniformly, $\delta$ is continuous. For every $n\in \N$, one proves, by induction, that $\rho_n\colon [0,1]^n\rightarrow [0,1]$ is monotone, and thus $\widetilde{\rho}_n\colon [0,1]^\N\rightarrow [0,1]$ is monotone, as well. Since $\delta$ is the pointwise limit of $\widetilde{\rho}_n$, $\delta$ is monotone, as well.
\endproof
The primitive operation symbols of $\MCinfty$ ($\delta$, $\land, \lor,\oplus, \odot$, and, for every $\lambda\in [0,1]$, the constant symbol $\lambda$) have a natural interpretation in $[0,1]$; indeed, each of them can be viewed as a function from a power of $[0,1]$ to $[0,1]$ itself. As we noticed in Remark~\ref{r:all operations in V are mon and cont} and in Lemma~\ref{l:delta-cont-mon}, they are monotone and continuous. In fact, the operations of $\MCinfty$ are \emph{all} monotone continuous functions from some power of $[0,1]$ to $[0,1]$. The following theorems make this statement precise.
\theorem \label{t:operations concide with monotone continuous}
	For each cardinal $\kappa$, the set of monotone continuous functions from $[0,1]^\kappa$ to $[0,1]$ coincides with the set of interpretations in $[0,1]$ of MC$_\infty$-terms of arity $\kappa$.
\endtheorem
\proof 
	Let $A$ be the set of functions $f\colon [0,1]^\kappa\rightarrow [0,1]$ for which there exist an MC$_\infty$-term (depending on $f$) of arity $\kappa$ whose interpretation in $[0,1]$ is $f$. Since the interpretation in $[0,1]$ of an MC$_\infty$-term is monotone and continuous by Remark~\ref{r:all operations in V are mon and cont} and Lemma~\ref{l:delta-cont-mon}, we have $A\subseteq \C\left([0,1]^\kappa\right)$. Moreover, $A$ contains, for each $i\in \kappa$, the projection $\pi_i\colon [0,1]^\kappa\rightarrow [0,1]$. Then, Theorem~\ref{t:StWe} applies, and we obtain that $A$ is dense in $\C\left([0,1]^\kappa\right)$. Furthermore, $A$ is an MC$_\infty$-algebra, and therefore it is Cauchy complete; thus $A=\C\left([0,1]^\kappa\right)$.
\endproof 

Roughly speaking, Theorem~\ref{t:operations concide with monotone continuous} says that the interpretation in $[0,1]$ is a surjective operator from the class of equivalence classes of MC$_\infty$-terms (where the equivalence relation is defined in the standard manner by identifying two terms if their interpretation in each algebra of the variety coincides) to the class of monotone continuous functions from some power of $[0,1]$ to $[0,1]$ itself. One consequence of the Theorem~\ref{t:W=ISP[0,1]} below is that this operator is injective, too, and so the equivalence classes of MC$_\infty$-terms are in bijective correspondence with the monotone continuous functions from some power of $[0,1]$ to $[0,1]$ itself. To state the theorem, we recall the standard operators $\mathrm{I}$ (closure under isomorphisms), $\mathrm{S}$ (closure under subalgebras) and $\mathrm{P}$ (closure under products). Moreover, we denote simply with $[0,1]$ the canonical MC$_\infty$-algebra whose underlying set is the unit interval $[0,1]$. 
\theorem \label{t:W=ISP[0,1]}
	$$\MCinfty=\mathrm{ISP}([0,1]).$$
\endtheorem
\proof 
	The right-to-left inclusion $\supseteq$ is clear beacuse $\MCinfty$ is a variety containing $[0,1]$. For the converse inclusion, let $A\in\MCinfty$. Then $A$ is archimedean, and thus the MC-morphisms towards $[0,1]$ separate the points of $A$. Since $\mathrm{U}_{\MCinfty,\MC}$ is full, every MC-morphism from $A$ to $[0,1]$ is also an MC$_\infty$-morphism. Hence, there are enough MC$_\infty$-morphisms from $A$ to $[0,1]$ to separate the points of $A$, and so $A$ is isomorphic to a subalgebra of a power of $[0,1]$.
\endproof 
\theorem 
	Let $I$ be a set. The MC$_\infty$-algebra $\C\left([0,1]^I\right)$ is freely generated by the projections $\left(\pi_i\colon \C\left([0,1]^I\right)\rightarrow [0,1]\right)_{i\in I}$.
\endtheorem
\proof 
	By Theorem~\ref{t:W=ISP[0,1]}, the free algebra generated by the set $I$ is the set of functions $[0,1]^I\rightarrow [0,1]$ that are the evaluation of a term of arity $\lvert I\rvert$. By Theorem~\ref{t:operations concide with monotone continuous}, this set is precisely $\C([0,1]^I)$.
\endproof 

Let us recall, from \citep[Section 1]{Linton}, Linton's definition of equational theory, varietal theory, equational category and varietal category. An \emph{equational theory} is a product preserving covariant functor $T\colon \Set^\op\rightarrow \mathbb{T}$ from the dual of the category of sets to a category $\mathbb{T}$ whose class of objects is put by $T$ in one-one correspondence with the objects of $\Set^\op$. One may then identify each object $\widetilde{n}$ of $\mathbb{T}$ with the set $n\in \Set$ such that $T(n)=\widetilde{n}$. The idea behind this definition is that the morphisms in $\mathbb{T}$ from $T(n)$ to $T(m)$ are the $\lvert m\rvert$-tuples of equivalence classes of terms of arity $\lvert n\rvert$ of a certain variety of algebras---where the equivalence relation is defined in the standard manner by identifying two terms if their interpretation in each algebra of the variety coincides---and the composition of morphisms is just the composition of terms (modulo the equivalence relation). From the category $\Set^\mathbb{T}$ of set valued functors on $\mathbb{T}$, we single out the full subcategory $\Set^T$ whose objects are the functors $X\colon \mathbb{T}\to \Set$ such that the composite $XT\colon \Set^\op\rightarrow \Set$ preserves products. One such functor $X$ is called a $T$-algebra. Any category equivalent to the category $\Set^T$ is called an \emph{equational category}.  Evaluation at the object $T(1)\in \mathbb{T}$ provides a faithful functor $U_T\colon \Set^T\rightarrow \Set$, the underlying set functor for $T$-algebras. The equational theory $T$ is called \emph{varietal} if the category $\mathbb{T}$ is locally small, and in this case any category equivalent to $\Set^T$ is said to be a \emph{varietal category}. 

Linton's setting generalizes Lawvere's perspective for finitary algebras \citep{Lawvere} to the infinitary ones (see \citep{Slo}). In fact, every variety of algebras $\mathbf{V}$ is a varietal category: $\mathbf{V}$ is equivalent to the category of $T$-algebras, where $\mathbb{T}$ is the opposite of the category of free algebras with homomorphisms, and $T\colon \Set^\op\rightarrow \mathbb{T}$ maps a set $I$ to the free algebra $\mathrm{Free}(I)$ over $I$. Note that the set of homomorphisms from $\mathrm{Free}(n)$ to $\mathrm{Free}(m)$ is in bijection with the set of $\lvert n\rvert$-tuples of equivalence classes of terms of arity $\lvert m\rvert$, where the equivalence relation is defined in the standard manner by identifying two terms if their interpretation in each algebra of the variety coincides.

\remark
The results in this section show that $\MCinfty$ is the category of algebras of the varietal theory $T\colon \Set \to \mathbb{T}$, where, for each set $I$, $T(I)= [0,1]^I$, and the morphisms from $T(I)$ to $T(J)$ are the monotone continuous maps from $[0,1]^I$ to $[0,1]^J$. The fact that the concrete varietal category $\Set^T$, i.e.\ $\MCinfty$, has a class of primitive operations of countable arity is equivalent to the fact that every continuous map from a power of $[0,1]$ to $[0,1]$ depends on at most countably many coordinates \citep[Theorem 1]{Mibu}. However, in this paper we do not settle the question whether $\PosComp^\op$ is equivalent or not to a variety of \emph{finitary} algebras; what we can say is that the functor $\hom(-,[0,1])\colon \PosComp\to \Set$ cannot be naturally isomorphic to the forgeftul functor of a variety of finitary algebras, because the function $\delta$ fails to be dependent on at most finitely many coordinates.
\endremark

\section{Conclusions}

$\C$ and $\Max$ establish a dual adjunction between $\PreT$ and $\MC$, induced by the dualizing object $[0,1]$. The fixed objects of this adjunction are precisely the objects in the images of the two functors: the fixed objects in $\PreT$ are the compact ordered spaces, while the fixed objects in $\MC$ are the archimedean Cauchy complete algebras. The forgetful functor from the variety $\MCinfty$ to the full subcategory of $\MC$ of archimedean Cauchy complete algebras is an isomorphism of categories. Therefore, $\C$ and $\Max$ restrict to a dual equivalence between $\PosComp$ and $\MCinfty$, induced by the dualizing object $[0,1]$. The main result is Theorem~\ref{t:MAIN-PosComp-is-variety}, i.e., the following.
\begin{quote}
	The category $\PosComp^\op$ is equivalent to a variety of algebras.
\end{quote}
The additional results are the description of the variety by means of operations and equational axioms, the description of the dual equivalence via the dualizing object $[0,1]$, and the extension of the duality to a wider dual adjunction between the category $\PreT$ and the finitary variety $\MC$.



\begin{thebibliography}{13}
	\providecommand{\natexlab}[1]{#1}
	\providecommand{\url}[1]{\texttt{#1}}
	\expandafter\ifx\csname urlstyle\endcsname\relax
	\providecommand{\doi}[1]{doi: #1}\else
	\providecommand{\doi}{doi: \begingroup \urlstyle{rm}\Url}\fi
	
	\bibitem[Ad\'{a}mek et~al.(2006)Ad\'{a}mek, Herrlich, and Strecker]{Joyofcats}
	J.~Ad\'{a}mek, H.~Herrlich, and G.~E. Strecker.
	\newblock Abstract and concrete categories: the joy of cats.
	\newblock \emph{Repr. Theory Appl. Categ.}, \penalty0 (17):\penalty0 1--507,
	2006.
	\newblock URL \url{http://www.tac.mta.ca/tac/reprints/articles/17/tr17.pdf}.
	\newblock Reprint of the 1990 original [Wiley, New York].
	
	\bibitem[Duskin(1969)]{Duskin}
	J.~Duskin.
	\newblock Variations on {B}eck's tripleability criterion.
	\newblock In \emph{Reports of the {M}idwest {C}ategory {S}eminar, {III}}, pages
	74--129. Springer, Berlin, 1969.
	
	\bibitem[Gierz et~al.(2003)Gierz, Hofmann, Keimel, Lawson, Mislove, and
	Scott]{ContLattDom}
	G.~Gierz, K.H. Hofmann, K.~Keimel, J.D. Lawson, M.~Mislove, and D.~S. Scott.
	\newblock \emph{Continuous lattices and domains}, volume~93 of
	\emph{Encyclopedia of Mathematics and its Applications}.
	\newblock Cambridge University Press, Cambridge, 2003.
	\newblock ISBN 0-521-80338-1.
	\newblock \doi{10.1017/CBO9780511542725}.
	\newblock URL \url{https://doi.org/10.1017/CBO9780511542725}.
	
	\bibitem[Hofmann and Nora(2018)]{HofmannLong}
	D.~Hofmann and P.~Nora.
	\newblock Enriched {S}tone-type dualities.
	\newblock \emph{Advances in Mathematics}, 330, 2018.
	\newblock \doi{10.1016/j.aim.2018.03.010}.
	
	\bibitem[Hofmann et~al.(2018)Hofmann, Neves, and Nora]{HofmannShort}
	D.~Hofmann, R.~Neves, and P.~Nora.
	\newblock Generating the algebraic theory of ${C}({X})$: the case of partially
	ordered compact spaces.
	\newblock \emph{Theory and Applications of Categories}, 33\penalty0 (12), 2018.
	
	\bibitem[Isbell(1982)]{Isbell}
	J.~Isbell.
	\newblock Generating the algebraic theory of {$C(X)$}.
	\newblock \emph{Algebra Universalis}, 15\penalty0 (2):\penalty0 153--155, 1982.
	\newblock ISSN 0002-5240.
	\newblock \doi{10.1007/BF02483718}.
	\newblock URL \url{https://doi.org/10.1007/BF02483718}.
	
	\bibitem[Lawvere(1963)]{Lawvere}
	F.~W. Lawvere.
	\newblock Functorial semantics of algebraic theories.
	\newblock \emph{Proc. Nat. Acad. Sci. U.S.A.}, 50:\penalty0 869--872, 1963.
	\newblock ISSN 0027-8424.
	\newblock \doi{10.1073/pnas.50.5.869}.
	\newblock URL
	\url{https://doi-org.pros.lib.unimi.it:2050/10.1073/pnas.50.5.869}.
	
	\bibitem[Linton(1966)]{Linton}
	F.~E.~J. Linton.
	\newblock Some aspects of equational categories.
	\newblock In \emph{Proc. {C}onf. {C}ategorical {A}lgebra ({L}a {J}olla,
		{C}alif., 1965)}, pages 84--94. Springer, New York, 1966.
	
	\bibitem[Marra and Reggio(2017)]{MarraReggio}
	V.~Marra and L.~Reggio.
	\newblock Stone duality above dimension zero: axiomatising the algebraic theory
	of {${\rm C}(X)$}.
	\newblock \emph{Adv. Math.}, 307:\penalty0 253--287, 2017.
	\newblock ISSN 0001-8708.
	\newblock \doi{10.1016/j.aim.2016.11.012}.
	\newblock URL \url{https://doi.org/10.1016/j.aim.2016.11.012}.
	
	\bibitem[Mibu(1944)]{Mibu}
	Y.~Mibu.
	\newblock On {B}aire functions on infinite product spaces.
	\newblock \emph{Proc. Imp. Acad. Tokyo}, 20:\penalty0 661--663, 1944.
	\newblock ISSN 0369-9846.
	\newblock URL
	\url{http://projecteuclid.org.pros.lib.unimi.it/euclid.pja/1195572745}.
	
	\bibitem[Nachbin(1965)]{Nachbin}
	L.~Nachbin.
	\newblock \emph{Topology and order}.
	\newblock Translated from the Portuguese by Lulu Bechtolsheim. Van Nostrand
	Mathematical Studies, No. 4. D. Van Nostrand Co., Inc., Princeton,
	N.J.-Toronto, Ont.-London, 1965.
	
	\bibitem[Porst and Tholen(1991)]{Tholen}
	H.E. Porst and W.~Tholen.
	\newblock Concrete dualities.
	\newblock In \emph{Category theory at work ({B}remen, 1990)}, volume~18 of
	\emph{Res. Exp. Math.}, pages 111--136. Heldermann, Berlin, 1991.
	
	\bibitem[S{\l}omi\'{n}ski(1959)]{Slo}
	J.~S{\l}omi\'{n}ski.
	\newblock The theory of abstract algebras with infinitary operations.
	\newblock \emph{Rozprawy Mat.}, 18:\penalty0 67 pp. (1959), 1959.
	\newblock ISSN 0860-2581.
	
\end{thebibliography}

\end{document}